\documentclass[11pt]{article} 
 
\usepackage{latexsym} 
\usepackage{amsfonts} 
\usepackage{amsmath} 
\usepackage{amsthm} 
\usepackage{mathrsfs} 
\usepackage{dsfont} 
\usepackage{bbold} 
\usepackage[english]{babel} 
\usepackage{caption} 
\usepackage{epsfig} 
\usepackage{float} 
\usepackage{subfigure} 
\usepackage{psfrag} 
\usepackage{graphicx} 
\usepackage{epsfig} 
\usepackage{amsbsy} 
\usepackage{hyperref}
 
\DeclareMathOperator{\Val}{\matV}

\newtheorem{theorem}{Theorem} 
\newtheorem*{prop*}{Theorem}

\newtheorem{defi}[theorem]{Definition} 
\newtheorem{lemma}[theorem]{Lemma} 
 
\newtheorem{rmk}[theorem]{Remark}

\newtheorem{cond}{Condition} 
 
\newcommand{\zerarcounters}{\setcounter{equation}{0}\setcounter{theorem}{0}} 
 
\newcommand{\ZZZ}{\mathds{Z}} 
\newcommand{\CCC}{\mathds{C}} 
\newcommand{\NNN}{\mathds{N}} 
 
\newcommand{\RRR}{\mathds{R}} 
\newcommand{\TTT}{\mathds{T}} 
\newcommand{\uno}{\mathds{1}} 
 
\newcommand{\BB}{{\mathcal B}} 
 
\newcommand{\DD}{{\mathcal D}} 
 
\newcommand{\calF}{{\mathcal F}} 
\newcommand{\calG}{{\mathcal G}}

\newcommand{\MM}{{\mathcal M}}

\newcommand{\calP}{{\mathcal P}} 
 
\newcommand{\RR}{{\mathcal R}} 
\newcommand{\SSSS}{{\mathcal S}}

\newcommand{\gotn}{{\mathfrak n}}

\newcommand{\gotB}{{\mathfrak B}} 
 
\newcommand{\gotD}{{\mathfrak D}} 
 
\newcommand{\gotF}{{\mathfrak F}}

\newcommand{\gotM}{{\mathfrak M}} 
\newcommand{\gotN}{{\mathfrak N}}

\newcommand{\gotR}{{\mathfrak R}} 
\newcommand{\gotS}{{\mathfrak S}}

\newcommand{\matD}{{\mathscr D}}

\newcommand{\matL}{{\mathscr L}}

\newcommand{\matR}{{\mathscr R}}

\newcommand{\matV}{{\mathscr V}}

\newcommand{\ol}{\overline} 
 
\newcommand{\Fullbox}{{\rule{2.0mm}{2.0mm}}} 
 
\newcommand{\EP}{\hfill\Fullbox\vspace{0.2cm}} 
\newcommand{\prova}{\noindent{\it Proof. }} 
\newcommand{\io}{\infty} 
\newcommand{\e}{\varepsilon} 
\newcommand{\al}{\alpha} 
 
\newcommand{\be}{\beta}

\newcommand{\x}{\xi}

\newcommand{\g}{\gamma}

\newcommand{\ze}{\zeta}

\newcommand{\s}{\sigma} 
 
\newcommand{\del}{\partial}

\newcommand{\oo}{\boldsymbol{\omega}} 
 
\newcommand{\aaa}{\boldsymbol{\alpha}} 
 
\newcommand{\ttheta}{\boldsymbol{\theta}} 
\newcommand{\hhh}{\boldsymbol{\eta}} 
\newcommand{\nn}{\boldsymbol{\nu}} 

\newcommand{\pps}{\boldsymbol{\psi}} 
\newcommand{\vzero}{\boldsymbol{0}}

\newcommand{\ttt}{\boldsymbol{t}}

\newcommand{\aaaa}{\boldsymbol{a}}

\newcommand{\FF}{\boldsymbol{F}}
\newcommand{\GG}{\boldsymbol{G}}
 
\newcommand{\AAA}{\boldsymbol{A}} 
\newcommand{\VVal}{\pmb{\Val}}
\newcommand{\TVal}{\overline{\Val}}
\newcommand{\TTVal}{\overline{\VVal}}
 
\newcommand{\ii}{{\rm i}}

\def\ins#1#2#3{\vbox to0pt{\kern-#2 \hbox{\kern#1 #3}\vss}\nointerlineskip} 
 
\begin{document}
 
\title{\bf Lower-dimensional invariant tori for perturbations\\
of a class of non-convex Hamiltonian functions} 
 
\author 
{\bf Livia Corsi$^{1}$, Roberto Feola$^{2}$ and Guido Gentile$^{2}$
\vspace{2mm} 
\\ \small
$^{1}$ Dipartimento di Matematica, Universit\`a di
Napoli ``Federico II'', Napoli, I-80126, Italy
\\ \small 
$^{2}$ Dipartimento di Matematica, Universit\`a di Roma Tre, Roma,
I-00146, Italy
\\ \small 
E-mail:  livia.corsi@unina.it, roberto\_feola@hotmail.com,
gentile@mat.uniroma3.it}

\date{} 
 
\maketitle 
 
\begin{abstract} 
We consider a class of quasi-integrable Hamiltonian systems obtained
by adding to a non-convex Hamiltonian function of an integrable system
a perturbation depending only on the angle variables.
We focus on a resonant maximal torus of the unperturbed system, 
foliated into a family of lower-dimensional tori of codimension 1,
invariant under a quasi-periodic flow with rotation vector
satisfying some mild Diophantine condition.
We show that at least one lower-dimensional torus with that
rotation vector always exists also for the perturbed system.
The proof is based on multiscale analysis and resummation procedures of divergent series.
A crucial role is played by suitable symmetries and cancellations,
ultimately due to the Hamiltonian structure of the system.
\end{abstract} 
  
\zerarcounters 
\section{Introduction} 
\label{sec:1} 

Consider the Hamiltonian dynamical system described, in action-angle variables,
by the Hamiltonian function
\begin{equation}\label{eq:1.1}
H(\aaa,\be,\AAA,B)=-\frac{1}{2}\AAA\cdot\AAA + \frac{1}{2}B^{2}+
\e f(\aaa,\be),
\end{equation}
where $(\aaa,\be)\in\TTT^{d}\times\TTT$, $(\AAA,B)\in\RRR^{d}\times\RRR$,
$f:\TTT^{d+1}\to\RRR$ is real-analytic, $\e\in\RRR$ is a small parameter
(the \emph{perturbation parameter})
and $\cdot$ is the standard scalar product in $\RRR^{d}$.
The corresponding Hamilton equations can be written as closed equations
for the angle variables $(\aaa,\be)$,
\begin{equation}\label{eq:1.2}
\left\{
\begin{aligned}
\ddot{\aaa}&=\e\partial_{\aaa}f(\aaa,\be),\\
\ddot{\be}&=-\e\partial_{\be}f(\aaa,\be).
\end{aligned}
\right.
\end{equation}

For $\e=0$ all the solutions of \eqref{eq:1.2} are trivially of the form
$(\aaa(t),\be(t))=(\aaa_{0}-\AAA_{0}t,\be_{0}+B_{0}t)$ where $(\aaa_{0},\be_{0})$ 
and $(\AAA_{0},B_{0})$ are the initial phases and actions, respectively.
Fix the initial actions as $(\AAA_{0},B_{0})=(-\oo,0)$
with $\oo\in\RRR^{d}$ such that $\oo\cdot\nn\ne0$ for all $\nn\in\ZZZ^{d}_{*}:=
\ZZZ^{d}\setminus\{\vzero\}$. Then the solutions
of the unperturbed system lie on a $(d+1)$-dimensional invariant torus
foliated into $d$-dimensional invariant tori parametrized by $\be_{0}$.

For $\e\neq0$ we say that the system has an invariant $d$-dimensional torus with
frequency $\oo$ if there is an invariant manifold for (\ref{eq:1.2})
where the motion is conjugated to a rotation with frequency vector $\oo$ on $\TTT^{d}$,
more precisely if there exists $\be_{0}\in\TTT$ and two  analytic functions
$\aaa_{\e}:\TTT^{d}\to\TTT^{d}$
and $\be_{\e}:\TTT^{d}\to\TTT$ such that
$\aaa_{0}(\pps)=\vzero$ and $\be_{0}(\pps)=0$,
the submanifold $\gotM$ of the form $\aaa=\pps+ \aaa_{\e}(\pps)$
and $\be=\be_{0} + \be_{\e}(\pps)$ is invariant
for (\ref{eq:1.2}) and the flow on $\gotM$ is given by $\pps\to \pps+\oo t$.

For $\oo\in\RRR^{d}$ define the \emph{Bryuno function} as \cite{B}
\begin{equation} \nonumber 
\BB(\oo):=
\sum_{m\ge0}\frac{1}{2^{m}}\log\frac{1}{\al_{m}(\oo)},\qquad
\al_{m}(\oo):=\inf_{\substack{\nn\in\ZZZ^{d}\\0<|\nn|\le 2^{m}}}|\oo\cdot\nn|.
\end{equation}

We shall prove the following result.

\begin{theorem}\label{thm:1}
For any $\oo\in\RRR^{d}$ such that $\BB(\oo)<\io$,
there exists $\e_{0}>0$ such that for any $|\e|<\e_{0}$ the system
(\ref{eq:1.2}) admits at least one invariant $d$-dimensional torus  with frequency $\oo$.
\end{theorem}

\begin{rmk} \label{rmk:1.2}
\emph{
It will turn out from the proof that it may happen that
the system admits a whole $(d+1)$-dimensional torus foliated
into $d$-dimensional invariant tori. In such a highly
non-generic case, the solution is analytic
in both the initial data and in the perturbation parameter.
}
\end{rmk}

Theorem \ref{thm:1} can be seen as a particular case of the result
announced in \cite{PK}, where the problem of existence of $d$-dimensional tori
is considered for Hamiltonian systems with $(d+1)$-degrees of freedom
described by Hamiltonian funtions $H(\aaa,\be,I)=H_{0}(I)+\e H_{1}(\aaa,\be,I)$,
where $I=(\AAA,B) \in \gotD_{1}$, with $\gotD_{1}$ a neighbourhood
of zero in $\RRR^{d+1}$, and the functions $H_{0},H_{1}$ are real-analytic in
all their arguments and $2\pi$-periodic in $\aaa,\be$ and
(modulo a canonical transformation) satisfy the following conditions:

\begin{enumerate}

\item \label{diof} $\partial_{\AAA}H_{0}(\vzero,0)=\oo$, with $\oo$ such that
$|\oo\cdot\nn| \ge c\,|\nn|^{-\tau}$ for $c>0$, $\tau>d-1$ and all $\nn\in\ZZZ^{d}_{*}$;

\item $\det \partial^{2}_{I} H_{0}(0) \neq 0$ and
$\partial_{B}^{2} H_{0}(0)=1$;

\item the function $H_{0}(\AAA,B)-\oo\cdot\AAA$ has a saddle point of
signature 1 in zero, that is 
$$ \hhh \cdot S_{0} \hhh < | \ttt_{0} \cdot \hhh|^{2}  \qquad
\forall\hhh\in\RRR^{d}\setminus\{\vzero\} ,
$$
where $S_{0}:=\partial^{2}_{\AAA} H_{0}(0)$ and $\ttt_{0}:=\partial_{\AAA}\partial_{B}H_{0}(0)$.

\end{enumerate}

\noindent 
Indeed the Hamiltonian function (\ref{eq:1.1}) satisfies the conditions above,
with $S_{0}=-\uno$ and $\ttt_{0}=\vzero$ (in fact the Bryuno condition $\BB(\oo)<\infty$
is weaker than the standard Diophantine condition in item \ref{diof} above).
Our method should apply to the more general case considered in \cite{PK}: we prefer to
focus on a particular class of systems, to avoid technical complications
and  put emphasis on the method, rather than the result itself --- already stated
in \cite{PK}. We shall show that the result can be credited
to the existence of remarkable symmetries of suitable quantities,
the so-called self-energies, that will be introduced in the proof.
In turn such symmetries are related to the Hamiltonian form of the equations of motion.

The construction envisaged below, as well as the method of \cite{PK}, 
does not allow us to obtain the existence of invariant $d$-dimensional tori
in the case of convex $H_{0}(I)$ treated by Cheng \cite{Ch}.
At the end we shall try to briefly illustrate where problems arise
when dealing with convex Hamiltonians. In particular we shall see
that the aforementioned symmetries are not sufficient in that case,
and other cancellation mechanisms should be looked for.

\zerarcounters 
\section{The formal expansion} 
\label{sec:2} 

Fix $\oo\in\RRR^{d}$ such that $\BB(\oo)<\io$. We look for a
quasi-periodic solution to \eqref{eq:1.2} of the form
$(\aaa(t),\be(t))=(\aaa_{0}+\oo t +\aaaa(t),\be_{0}+b(t))$, with
\begin{equation} \label{eq:2.1}
\aaaa(t)=\sum_{\nn\in\ZZZ^{d}_{*}}{\rm e}^{\ii\nn\cdot\oo t}\aaaa_{\nn},
\qquad
b(t)=\sum_{\nn\in\ZZZ^{d}_{*}}{\rm e}^{\ii\nn\cdot\oo t}b_{\nn},
\end{equation}
and $(\aaaa(t),b(t))\to(\vzero,0)$ as $\e\to0$,
so that in the Fourier space \eqref{eq:1.2} becomes
\begin{subequations}
\begin{align}
&(\oo\cdot\nn)^{2}\aaaa_{\nn}=-[\e \partial_{\aaa} f(\aaa,\be)]_{\nn},
\qquad \nn\ne\vzero,
\label{eq:2.2a} \\
&(\oo\cdot\nn)^{2}b_{\nn}=[\e\partial_{\be} f(\aaa,\be)]_{\nn},
\qquad \nn\ne\vzero,
\label{eq:2.2b} \\
&[\e \partial_{\aaa} f(\aaa,\be)]_{\vzero}=\vzero,
\label{eq:2.2c} \\
&[\e \partial_{\be} f(\aaa,\be)]_{\vzero}=0 ,
\label{eq:2.2d}
\end{align}
\label{eq:2.2}
\end{subequations}
\vskip-.5truecm \noindent
where
\begin{equation}\label{eq:2.3}
\begin{aligned}
&[\partial_{\aaa}f(\aaa,\be)]_{\nn}=\sum_{\substack{p\geq0 \\ q\geq0}}
\sum_{\substack{\nn_0+\ldots+\nn_{p+q}=\nn \\ \nn_{0}\in\ZZZ^d \\ 
\nn_i \in\ZZZ^{d}_{*}, i=1,\ldots,p+q}}
\frac{1}{p!q!}
{{(\ii\nn_0)^{p+1}}}{{\partial_{\be}^{q}
f_{\nn_0}(\aaa_0,\be_{0})}}\prod_{i=1}^{p}\aaaa_{\nn_i}
\prod_{j=p+1}^{p+q}b_{\nn_j},\\
&[\partial_{\be}f(\aaa,\be)]_{\nn}=\sum_{\substack{p\geq0 \\ q\geq0}}
\sum_{\substack{\nn_0+\ldots+\nn_{p+q}=\nn \\ 
\nn_{0}\in\ZZZ^d \\ \nn_i \in\ZZZ^{d}_{*}, i=1,\ldots,p+q}}
\frac{1}{p!q!}
{{(\ii\nn_0)^{p}}}{{\partial_{\be}^{q+1}
f_{\nn_0}(\aaa_0,\be_{0})}}\prod_{i=1}^{p}\aaaa_{\nn_i}\prod_{j=p+1}^{p+q}
b_{\nn_j},
\end{aligned}
\end{equation}
\vskip-.1truecm \noindent
and $f_{\nn}(\aaa_0,\be_0)={\rm e}^{\ii\nn\cdot\aaa_0}\hat{f}_{\nn}(\be_0)$,
where we denoted
$$
f(\aaa,\be)=\sum_{\nn\in\ZZZ^d}\hat{f}_{\nn}(\be){\rm e}^{\ii\nn\cdot\aaa}.
$$
Throughout the paper, the sums and the products over the empty set
have to be considered as $0$ and $1$, respectively.
Equations \eqref{eq:2.2a} and \eqref{eq:2.2b} are called the
\emph{range equations}, while \eqref{eq:2.2c} and \eqref{eq:2.2d}
are called the \emph{bifurcation equations}.

We start by writing formally
\begin{subequations}
\begin{align}
&\aaa(t)=\aaa(t;\e,\aaa_0,\be_0)=\aaa_0+\oo t+\sum_{k\geq1}\e^{k}
\sum_{\nn\in\ZZZ^{d}_{*}}{\rm e}^{\ii\nn\cdot\oo t}\aaaa_{\nn}^{(k)}(\aaa_0,\be_0),
\label{eq:2.4a} \\
&\be(t)=\be(t;\e,\aaa_0,\be_0)=\be_0+\sum_{k\geq1}\e^{k}
\sum_{\nn\in\ZZZ^{d}_{*}}{\rm e}^{\ii\nn\cdot\oo t}b_{\nn}^{(k)}(\aaa_0,\be_0).
\label{eq:2.4b}
\end{align}
\label{eq:2.4}
\end{subequations}
\vskip-.1truecm \noindent 
If we define recursively for $k\ge1$ and $\nn\ne\vzero$
\begin{equation} \nonumber \label{eq:2.5}
\aaaa_{\nn}^{(k)}=-\frac{1}{(\oo\cdot\nn)^2}
[\partial_{\aaa}f(\aaa,\be)]_{\nn}^{(k-1)}, \qquad 
b_{\nn}^{(k)}=\frac{1}{(\oo\cdot\nn)^2}
[\partial_{\be}f(\aaa,\be)]_{\nn}^{(k-1)},
\end{equation}
with
$[\del_{\aaa}f(\aaa,\be)]_{\nn}^{(0)}=\ii\nn\,f_{\nn}(\aaa_0,\be_{0})$,
$[\del_{\be}f(\aaa,\be)]_{\nn}^{(0)}=\del_{\be}f_{\nn}(\aaa_0,\be_0)$ and, for $k\geq1$,
\begin{equation} \nonumber 
\begin{aligned}
&[\del_{\aaa}f(\aaa,\be)]_{\nn}^{(k)}=\sum_{\substack{p\geq0 \\ q\geq0}}
\sum_{\substack{\nn_0+\ldots+\nn_{p+q}=\nn \\ \nn_{0}\in\ZZZ^d \\ 
\nn_i \in\ZZZ^{d}_{*}, i=1,\ldots,p+q}}
\frac{1}{p!q!}
{{(\ii\nn_0)^{p+1}}}{{\partial_{\be}^{q}
f_{\nn_0}(\aaa_0,\be_{0})}}
\sum_{\substack{k_{1}+\ldots+k_{p+q}=k\\k_{i}\ge1}}
\prod_{i=1}^{p}\aaaa_{\nn_i}^{(k_{i})}
\prod_{j=p+1}^{p+q}b_{\nn_j}^{(k_{j})},\\
&[\partial_{\be}f(\aaa,\be)]_{\nn}^{(k)}=\sum_{\substack{p\geq0 \\ q\geq0}}
\sum_{\substack{\nn_0+\ldots+\nn_{s+r}=\nn \\ 
\nn_{0}\in\ZZZ^d \\ \nn_i \in\ZZZ^{d}_{*}, i=1,\ldots,p+q}}
\frac{1}{p!q!}
{{(\ii\nn_0)^{p}}}{{\partial_{\be}^{q+1}
f_{\nn_0}(\aaa_0,\be_{0})}}
\sum_{\substack{k_{1}+\ldots+k_{p+q}=k\\k_{i}\ge1}}
\prod_{i=1}^{p}\aaaa_{\nn_i}^{(k_{i})}\prod_{j=p+1}^{p+q}
b_{\nn_j}^{(k_{j})},
\end{aligned}
\end{equation}
then the series \eqref{eq:2.4} turn out to be a formal solution of
the range equations for any values of the parameters $\aaa_{0}$ and $\be_{0}$.

Unfortunately in general we are not able to prove the convergence of the
series \eqref{eq:2.4} and moreover we also have to solve the
bifurcation equations. As we shall see the two problems are somehow related.

\zerarcounters 
\section{Conditions of convergence for the formal expansion} 
\label{sec:3} 

In this section we shall see how to represent graphically the formal
solutions \eqref{eq:2.4}. We shall see that under suitable
(quite non-generic) hypotheses the two series converge.
However, in general, a resummation is needed to give the series a
meaning: this will be discussed in Section \ref{sec:4}.

\subsection{Diagrammatic rules}
\label{sec:3.1}

Our aim is to represent the formal series as a ``sum
over trees'', so first of all we need some definitions.
(We closely follow \cite{CG2,CG3}, with obvious adaptations).

A graph is a set of points and lines connecting them. 
A \emph{rooted tree} $\theta$ is a graph with no cycle, 
such that all the lines are oriented toward a single point 
(\emph{root}) which has only one incident line $\ell_{\theta}$ 
(\emph{root line});  we will omit the adjective ``rooted" in the following.
All the points in a tree except the root are called \emph{nodes}. 
The orientation of the lines in a tree induces a partial ordering 
relation ($\preceq$) between the nodes and the lines: we can 
imagine that each line carries an arrow pointing toward the root. 
Given two nodes $v$ and $w$, 
we shall write $w \prec v$ every time $v$ is along the path 
(of lines) which connects $w$ to the root. 

We denote by $N(\theta)$ and $L(\theta)$ the sets of nodes and 
lines in $\theta$, respectively. 
Since a line $\ell\in L(\theta)$ is uniquely identified 
by the node $v$ which it leaves, we may write $\ell = \ell_{v}$. 
We write $\ell_{w} \prec \ell_{v}$ if $w\prec v$, and $w\prec\ell=\ell_{v}$
if $w\preceq v$; 
if $\ell$ and $\ell'$ are two comparable lines, i.e. 
$\ell' \prec \ell$, we denote by $\calP(\ell,\ell')$ the 
(unique) path of lines connecting $\ell'$ to $\ell$, with $\ell$ and 
$\ell'$ not included (in particular $\calP(\ell,\ell')=\emptyset$ 
if $\ell'$ enters the node $\ell$ exits). 

With each node $v\in N(\theta)$ we associate a \emph{mode} label 
$\nn_{v}\in \ZZZ^{d}$ and a \emph{component} label
$h_{v}\in\{\al_{1},\ldots,\al_{d},\be\}$, and we denote by $s_{v}$ the
number of lines 
entering $v$. With each line $\ell=\ell_{v}$ we associate a
component label $h_{\ell_{v}}=h_{v}$ and 
a \emph{momentum} $\nn_{\ell}\in \ZZZ^{d}_{*}$, 
except for the root line which can have either zero momentum or not, 
i.e. $\nn_{\ell_{\theta}}\in\ZZZ^{d}$. For any node $v\in N(\theta)$
we denote by $p_{j,v},q_{v}$ the number of lines entering $v$ with component
$\al_{j}$ and $\be$, respectively, and set
$p_{v}=p_{1,v}+\ldots+p_{d,v}$; of course $s_{v}=p_{v}+q_{v}$.
Finally, we associate with each line $\ell$ also a \emph{scale label} 
such that $n_{\ell}=-1$ if $\nn_{\ell}=\vzero$, while 
$n_{\ell}\in\ZZZ_{+}$ if $\nn_{\ell}\neq\vzero$
(so far there is no relation between non-zero momenta and
scale labels: a constraint will appear shortly).
Note that one can have $n_{\ell}=-1$ only if $\ell$ is the root line
of $\theta$. 
We force the following \emph{conservation law} 
\begin{equation}\label{eq:3.1} 
\nn_{\ell}=\sum_{\substack{ w \in N(\theta) \\ w\prec \ell}}\nn_{w}. 
\end{equation} 

We shall call trees tout court the trees with labels, and we shall use
the term \emph{unlabelled tree} for the trees without labels. 
We shall say that two trees are \emph{equivalent} if they can be 
transformed into each other by continuously deforming the lines in 
such a way that these do not cross each other and also labels match. 
This provides an equivalence relation on the set of the trees. From
now on we shall call trees such equivalence classes. 

Given a tree $\theta$ we call \emph{order} of $\theta$ the 
number $k(\theta)=|N(\theta)|=|L(\theta)|$ (for any finite set $S$ 
we denote by $|S|$ its cardinality), \emph{total momentum} of 
$\theta$ the momentum associated with $\ell_{\theta}$
and \emph{total component} of $\theta$ the component associated with $\ell_{\theta}$. 
We shall denote by $\Theta_{k,\nn,h}$ the set of  trees 
with order $k$, total momentum $\nn$ and total component $h$.
A subset $T\subset\theta$ is a \emph{subgraph} of $\theta$
if it is formed by set of nodes $N(T)\subseteq
N(\theta)$ and lines $L(T)\subseteq L(\theta)$ connecting them
(possibly including the root line: in such a case we say that the
root is included in $T$) in such
a way that $N(T)\cup L(T)$ is connected.
If $T$ is a subgraph of $\theta$ we call \emph{order}
of $T$ the number $k(T)=|N(T)|$. We say that a line enters $T$ if it 
connects a node $v\notin N(T)$ to a node $w\in N(T)$, 
and we say  that a line exits $T$ if it connects a node $v\in N(T)$ 
to a node $w\notin N(T)$ or to the root (which is not included in $T$
in this case). Of course, if a line $\ell$ enters or exits $T$,
then $\ell\notin L(T)$. If $T$ is a subgraph 
of $\theta$ with only one entering line $\ell'$ and one exiting line 
$\ell$, we set $\calP_{T}:=\calP(\ell,\ell')$. 

A \emph{cluster} $T$ on scale $n$ is a maximal subgraph 
of a tree $\theta$ such that all the lines have scales 
$n'\le n$ and there is at least a line with scale $n$. 
The lines entering the cluster $T$ and the line coming 
out from it (unique if existing at all) are called the 
\emph{external} lines of $T$.
 
A \emph{self-energy cluster} is a cluster $T$ such that 
(i) $T$ has only one entering line $\ell'_{T}$ and 
one exiting line $\ell_{T}$, (ii) $\nn_{\ell}\ne\nn_{\ell'_{T}}$ for all
$\ell\in\calP_{T}$, (iii) one has $\nn_{\ell_{T}}= 
\nn_{\ell'_{T}}$ and hence $\sum_{v\in N(T)}\nn_{v}=\vzero$. 

We shall say that a self-energy cluster is on 
scale $-1$, if $N(T)=\{v\}$, with of course $\nn_{v}=\vzero$ 
(so that $\calP_{T}=\emptyset$). 

\begin{rmk}\label{rmk:3.1} 
\emph{ 
Given a self-energy cluster $T$, the momenta of the lines in $\calP_{T}$ 
depend on $\nn_{\ell'_{T}}$ because of the conservation law (\ref{eq:3.1}). 
More precisely, for all $\ell\in\calP_{T}$ one has 
$\nn_{\ell}=\nn_{\ell}^{0}+\nn_{\ell'_{T}}$ with 
$\nn_{\ell}^{0}=\sum_{\substack{w\in N(T) , w\prec \ell}} \nn_{w}$,
while all the other labels in $T$ do not depend on $\nn_{\ell'_{T}}$.
} 
\end{rmk} 

We say that two self-energy clusters $T_{1},T_{2}$ have the same 
\emph{structure} if setting $\nn_{\ell'_{T_{1}}}=\nn_{\ell'_{T_{2}}}=\vzero$ one has 
$T_{1}=T_{2}$. Of course this provides an equivalence relation on the 
set of all self-energy clusters. From now on we shall call 
self-energy clusters tout court such equivalence 
classes and we shall denote by $\gotS^{k}_{n,u,e}$ the set of
self-energy clusters with order $k$, scale $n$ and such that
$h_{\ell'_{T}}=e$ and $h_{\ell_{T}}=u$, with $e,u\in\{\al_{1},\ldots,\al_{d},
\be\}$.

Given any tree $\theta\in\Theta_{k,\nn,h}$ we associate with each node
$v\in N(\theta)$ a \emph{node factor}
\begin{equation}\label{eq:3.2}
\calF_{v}:=
\left\{
\begin{aligned}
& -\frac{1}{p_{v}!q_{v}!}{{(\ii\nn_{v})^{p_v+1}}}
{{\del_{\be}^{q_v}f_{\nn_{v}}(\aaa_0,\be_0)}}, \qquad 
h_v=\al_j, j=1,\ldots, d,\\
& \\
&\frac{1}{p_{v}!q_{v}!}  {{(\ii\nn_{v})^{p_{v}}}}
{{\del_{\be}^{q_v+1}f_{\nn_{v}}(\aaa_0,\be_0)}}, \qquad 
h_v=\be
\end{aligned}
\right.
\end{equation}
which is a tensor of rank $s_{v}+1$. We associate with each line
$\ell\in L(\theta)$ a \emph{propagator} defined as follows.
Let us introduce the sequences $\{m_{n},p_{n}\}_{n \ge 0}$, with $m_{0}=0$ 
and, for all $n\ge 0$, $m_{n+1}=m_{n}+p_{n}+1$,  where
$p_{n}:=\max\{q\in\ZZZ_{+}\,:\,\al_{m_{n}}(\oo)<2\al_{m_{n}+q}(\oo)\}$. Then the 
subsequence $\{\al_{m_{n}}(\oo)\}_{n\ge 0}$ of $\{\al_{m}(\oo)\}_{m\ge0}$ is
decreasing.  Let $\chi:\RRR\to\RRR$ be a $C^{\io}$ even function,
non-increasing for $x\ge0$, such that
\begin{equation}\label{eq:3.3} 
\chi(x)=\left\{ 
\begin{aligned} 
&1,\qquad |x| \le 1/2, \\ 
&0,\qquad |x| \ge 1. 
\end{aligned}\right. 
\end{equation} 
Set $\chi_{-1}(x)=1$ and $\chi_{n}(x)=\chi(8x/\al_{m_{n}}(\oo))$ for $n\ge0$. 
Set also $\psi(x)=1-\chi(x)$, $\psi_{n}(x)=\psi(8x/\al_{m_{n}}(\oo))$, 
and $\Psi_{n}(x)=\chi_{n-1}(x)\psi_{n}(x)$, for $n\ge 0$; see Figure 3.5 in \cite{CG3}.
Then we associate with each line a propagator
\begin{equation}\label{eq:3.4}
{\calG}_{\ell}:=\left\{
\begin{aligned}
&\frac{\Psi_{n_{\ell}}(\oo\cdot\nn_{\ell})}{(\oo\cdot\nn_{\ell})^{2}},
& \qquad n_{\ell}\ge0, \hskip.3truecm  & \\
&1,
& \qquad n_{\ell}=-1.
\end{aligned}
\right.
\end{equation}

Given any subgraph $S$ of any tree $\theta$ we define the
\emph{value} of $S$ as
\begin{equation}\label{eq:3.5}
\Val(S)=\left(\prod_{v\in N(S)}\calF_{v}\right)\left(
\prod_{\ell\in L(S)} \calG_{\ell}\right).
\end{equation}

Set $\Theta_{k,\nn,\aaa}:=\Theta_{k,\nn,\al_1}\times\ldots\times
\Theta_{k,\nn,\al_d}$ and for any
$\ttheta=(\theta_1,\ldots,\theta_d)\in\Theta_{k,\nn,\aaa}$ define 
$\VVal(\ttheta):=\left(\Val(\theta_1),\ldots,\Val(\theta_d)\right)$,
so that one has
\begin{subequations}
\begin{align}
&\aaaa_{\nn}^{(k)}= \!\!\!\!\!\!
\sum_{\ttheta\in\Theta_{k,\nn,\aaa}} \!\!\! \VVal(\ttheta),
\qquad b_{\nn}^{(k)}=\!\!\!\!\!\!
\sum_{\theta\in\Theta_{k,\nn,\be}} \!\!\! \Val(\theta),
\qquad  \nn\ne\vzero,
\label{eq:3.6a} \\
&[-\partial_{\aaa}f(\aaa,\be)]^{(k)}_{\vzero}=\!\!\!\!\!\!
\sum_{\ttheta\in\Theta_{k+1,\vzero\,\aaa}} \!\!\! \VVal(\ttheta),
\qquad  [\partial_{\be}f(\aaa,\be)]^{(k)}_{\vzero}=\!\!\!\!\!\!
\sum_{\theta\in\Theta_{k+1,\vzero,\be}} \!\!\! \Val(\theta),
\label{eq:3.6b}
\end{align}
\label{eq:3.6}
\end{subequations}
\vskip-.1truecm \noindent  
as is easy to check. In particular the quantities in \eqref{eq:3.6} are well
defined for any (fixed) $k\ge 1$ (see Appendix H in \cite{CG2}).

\begin{rmk}\label{rmk:3.2}
\emph{
Given a subgraph $S$ of any tree $\theta$ such that $\Val(S)\neq0$, for any line
$\ell\in L(S)$ (except possibly the root line of $\theta$) one has $\Psi_{n_{\ell}}
(\oo\cdot\nn_{\ell})\neq0$, so that
\begin{equation} \nonumber 
\frac{\al_{m_{n_\ell}}(\oo)}{16}\leq|\oo\cdot\nn_{\ell}|
\leq\frac{\al_{m_{n_{\ell}-1}}(\oo)}{8} <
\frac{\al_{m_{n_{\ell}-1}+p_{n_{\ell}-1}}(\oo)}{4}  = 
\frac{\al_{m_{n_{\ell}}-1}(\oo)}{4} ,
\end{equation}
where $\al_{m_{-1}}(\oo)$ has to be interpreted as $+\infty$,
and hence, by definition of $\al_{m}(\oo)$, one has $|\nn_{\ell}|> 
2^{m_{n_{\ell}}-1}$.
Moreover, by definition of $\{\al_{m_n}(\oo)\}_{n\geq0}$, the number of scales 
which can be associated with a line $\ell$ in such way that the propagator 
does not vanish is at most 2.
}
\end{rmk}

\subsection{Dimensional bounds}\label{sec:3.2}

For any subgraph $S$ of any tree 
$\theta$ call $\gotN_{n}(S)$ the number of lines on scale
$\geq n$ in $S$, and set
\begin{equation}\label{eq:3.7}
K(S):=\sum_{v\in N(S)}|\nn_{v}|.
\end{equation}

We shall say that a line $\ell$ is \emph{resonant} if it exits a self-energy cluster,
otherwise $\ell$ is \emph{non-resonant}. 
For any line $\ell \in \theta$ define the \emph{minimum scale} of $\ell$ as
\begin{equation*} 
\ze_{\ell} := \min\{ n \in \ZZZ_{+} : \Psi_{n} (\oo\cdot\nn_{\ell}) \neq 0 \} .
\end{equation*}
Given any subgraph $S$ of any tree $\theta$, we denote by
$\gotN^{\bullet}_{n}(S)$ the number of non-resonant lines
$\ell\in L(S)$ such that $\ze_{\ell} \ge n$. By definition, if $\Val(S)\neq0$,
for each line $\ell\in L(S)$ either $n_{\ell}=\ze_{\ell}$ or $n_{\ell}=\ze_{\ell}+1$.
We have the following results.

\begin{lemma}\label{lem:3.3}
For all $h\in\{\al_{1},\ldots,\al_{d},\be\}$, $\nn\in\ZZZ^{d}$,
$k\ge1$ and for any
$\theta\in\Theta_{k,\nn,h}$ with $\Val(\theta)\ne 0$, one has
$\gotN_{n}^{\bullet}(\theta)\le 2^{-(m_{n}-2)}K(\theta)$ for all $n\ge0$.
\end{lemma}

\begin{lemma}\label{lem:3.4}
For all $e,u\in\{\al_{1},\ldots,\al_{d},\be\}$, $n \ge 0$, $k\ge1$ and for any $T\in\gotS^{k}_{n,u,e}$
with $\Val(T)\ne 0$, one has $K(T)>2^{m_{n}-1}$ and
$\gotN_{p}^{\bullet}(T)\le 2^{-(m_{p}-3)}K(T)$ for all $0\le p\le n$.
\end{lemma}

The proofs of the two results above can be easily adapted from the proofs of Lemmas
6.4 and 6.5 in \cite{CG3}, respectively (and the same notations have been used),
notwithstanding the slightly different definition of resonant lines and the fact
that here the lines different from the root line can have only scale $\ge0$.

\begin{lemma}\label{lem:3.5}
For any tree $\theta\in\Theta_{k,\nn,h}$ and any self-energy cluster
$T\in\gotS^{k}_{n,u,e}$ denote by $L_{NR}(\theta)$ and $L_{NR}(T)$ the sets
of non-resonant lines in $\theta$ and $T$, respectively, and set
\begin{equation}\nonumber
\Val_{NR}(\theta):=\Biggl(\prod_{v\in N(\theta)}\calF_{v}\Biggr)
\Biggl(\prod_{\ell\in L_{NR}(\theta)}\calG_{\ell}\Biggr),\qquad
\Val_{NR}(T):=\Biggl(\prod_{v\in N(T)}\calF_{v}\Biggr)
\Biggl(\prod_{\ell\in L_{NR}(T)}\calG_{\ell}\Biggr) .
\end{equation}
Then
\begin{equation} \label{eq:3.8}
|\Val_{NR}(\theta)|
\le C_{1}^{k}{\rm e}^{-\xi|\nn|/2}, \qquad
|\Val_{NR}(T)|\le C_{2}^{k}{\rm e}^{-\xi K(T)/2},
\end{equation}
for some positive constants $C_{1}$ and $C_{2}$.
\end{lemma}

\prova
We prove only the first bound in (\ref{eq:3.8}) since the proof of
the second one proceeds in the
same way, with $T$ playing the role of $\theta$. For any $n_{0}\ge0$ one has
\begin{equation}\nonumber
\prod_{\ell\in L_{NR}(\theta)}|\calG_{\ell}| \le
\left(\frac{16}{\al_{m_{n_{0}}}(\oo)}\right)^{2k} 
\!\!\!\!\! \prod_{n\ge n_{0}+1}
\left(\frac{16}{\al_{m_{n}}(\oo)}\right)^{2\gotN_{n}^{\bullet}(\theta)}
\le D(n_{0})^{2k}{\rm exp}(\x(n_{0})K(\theta)), 
\end{equation}
with 
$$ 
D(n_{0})=\frac{16}{\al_{m_{n_{0}}}(\oo)},\qquad 
\x(n_{0})=8\sum_{n\ge n_{0}+1}\frac{1}{2^{m_{n}}}\log
\frac{16}{\al_{m_{n}}(\oo)}. 
$$ 
Then, since $\BB(\oo)<\io$, one can choose $n_{0}$ such that 
$\x(n_{0})\le \x/2$, so that, since
$$
\prod_{v\in N(\theta)}|\calF_{v}|\le C_{0}^{k}{\rm e}^{-\xi K(\theta)},
$$
for some positive constant $C_{0}$, the bound  follows.\EP

If $T$ is a self-energy cluster, we can (and shall)
write $\Val(T)=\Val_{T}(\oo\cdot\nn_{\ell'_{T}})$ and
$\Val_{NR}(T)=\Val_{T,NR}(\oo\cdot\nn_{\ell'_{T}})$
to stress the dependence on $\nn_{\ell'_{T}}$ --- see Remark \ref{rmk:3.1}.

\begin{rmk}\label{rmk:3.6}
\emph{
Since the proofs of Lemmas \ref{lem:3.3} and \ref{lem:3.4} work under the
weaker condition
$$
\frac{\al_{m_{n_{\ell}}}(\oo)}{32}<|\oo\cdot\nn_{\ell}|<
\frac{\al_{m_{n_{\ell}-1}}(\oo)}{4}
$$
one can show that also
$\partial_{x}^{j}\Val_{T,NR}(\tau x)$ admits the same bound
as $\Val_{T,NR}(x)$ in (\ref{eq:3.8}) for $j=0,1,2$ and
$\tau\in[0,1]$, possibly with a different constant $C_{2}$.
}
\end{rmk}
 
What emerges from Lemma \ref{lem:3.5} is that, if we could
ignore the resonant lines, the convergence of the series (\ref{eq:2.4})
would immediately follow (for $\e$ small enough). On the contrary,
the presence of resonant lines may be a real obstruction for the convergence.
Suppose indeed that a resonant line $\ell$ exits a self-energy cluster $T$
on scale $n \ll n_{\ell}$. Then $T$ must contain at least one line
$\ell'$ on scale $n$ such that $ 2|\oo\cdot\nn_{\ell'}^{0}| \ge
|\oo\cdot\nn_{\ell'}| \ge |\oo\cdot\nn_{\ell'}^{0}|/2$ since $n\ll n_{\ell}$
(recall the definition of $\nn_{\ell'}^{0}$ in  Remark \ref{rmk:3.1})
and hence $|\nn_{\ell'}^{0}| \ge 2^{m_{n}-1}$
(reason as in Remark \ref{rmk:3.2} to bound $|\nn_{\ell'}^{0}|$
in terms of the scale $n_{\ell}$). Therefore we can extract a factor
${\rm e}^{- \xi 2^{m_{n}}/8}$ from the product of the node factors of the nodes in $T$:
however this is not enough to control the
propagator $\calG_{\ell}$ for which we only have the bound
$2^{8}/\al_{m_{n_{\ell}}}(\oo)^{2}$. Moreover in principle a tree can contain
a ``chain'' of self-energy clusters
and hence of resonant lines, which implies accumulation of small
divisors. Therefore one would need a ``gain factor'' proportional to
$(\oo\cdot\nn_{\ell})^{2}$ for each resonant line $\ell$ for
the power series (\ref{eq:2.4}) to converge.

\subsection{Symmetries}\label{sec:3.3}

For all $k\ge1$ define the \emph{self-energies}
\begin{equation}\label{eq:3.9}
\begin{aligned}
M^{(k)}_{u,e}(x,n) & := \!\!\!\! \sum_{T\in\gotS^{k}_{n,u,e}} \!\!\!\! \Val_{T}(x),
\qquad \MM^{(k)}_{u,e}(x,n):= \!\! \sum_{p=-1}^{n}M^{(k)}_{u,e}(x,p) , \\
\MM^{(k)}_{u,e}(x) & :=\lim_{n\to\io}\MM^{(k)}_{u,e}(x,n).
\end{aligned}
\end{equation}
\vskip-.3truecm \noindent
Here we shall exhibit the existence of suitable symmetries for the self-energy
clusters, i.e some remarkable identities between the quantities
$\MM^{(k)}_{u,e}(x,n)$ and $\MM^{(k)}_{u,e}(x)$ introduced in (\ref{eq:3.9}).
In turn such symmetries will allow us to obtain a gain factor proportional to
$(\oo\cdot\nn_{\ell})^{2}$ for ``some'' resonant line $\ell$: under suitable
assumptions (which we shall exploit later on) this will imply the convergence
of the power series \eqref{eq:2.4}.

\begin{lemma}\label{lem:3.7}
For all $k\ge1$ one has
\begin{equation} \nonumber 
\null \hskip1.truecm 
\begin{matrix}
\MM^{(k)}_{\al_{i},\al_{j}}(0)
= \partial_{\al_{0,j}}[-\partial_{\al_{i}}f(\aaa,\be)]^{(k)}_{\vzero} , & \qquad
\MM^{(k)}_{\al_{i},\be}(0) = \partial_{\be_{0}}[-\partial_{\al_{i}}f(\aaa,\be)]^{(k)}_{\vzero} , \\
\MM^{(k)}_{\be,\al_{j}}(0) =
\partial_{\al_{0,j}}[\partial_{\be}f(\aaa,\be)]^{(k)}_{\vzero} , & \qquad
\MM^{(k)}_{\be,\be}(0) = \partial_{\be_{0}}[\partial_{\be}f(\aaa,\be)]^{(k)}_{\vzero} .
\end{matrix}
\end{equation}
\vskip-.3truecm
\end{lemma}

\prova
First of all let us write, for $e_{0}=\al_{0,1},\ldots,\al_{0,d},\be_{0}$
and $u=\al_{1},\ldots,\al_{d},\be$,
\begin{equation}\label{eq:3.10}
\partial_{e_{0}}\left(\sum_{\theta\in\Theta_{k,\vzero,u}}\Val(\theta)\right)
=\sum_{\theta\in\Theta_{k,\vzero,u}} \sum_{v\in N(\theta)}
\partial_{e_{0}}\calF_{v}\left( \prod_{v'\in N(\theta)\setminus\{v\}}
\calF_{v'}\right) \left(\prod_{\ell\in L(\theta)}\calG_{\ell}\right)  ,
\end{equation}
where we have used the fact that $\Val(\theta)$
depends on $\aaa_{0},\be_{0}$ only through the node factors.
Each summand in the r.h.s. of (\ref{eq:3.10}) differs from
$\Val(\theta)$ because a further derivative (with respect to
$\al_{0,j}$ or $\be_{0}$) acts on the node factor of a node $v\in N(\theta)$.
This can be graphically represented as the 
same tree $\theta$, but with a further line $\ell'$ entering the node $v$;
such a line carries $\vzero$-momentum and has component $e=\al_{1},\ldots,
\al_{d},\be$ for $e_{0}=\al_{0,1},\ldots,\al_{0,d},\be_{0}$, respectively, and
hence it is a contribution to $\MM^{(k)}_{u,e}(0)$. On the other hand
it is easy to realise that each contribution to $\MM^{(k)}_{u,e}(0)$ is of
the form above. Therefore the assertion follows.
\EP

\begin{lemma}\label{lem:3.8}
For all $k\ge1$ one has
\begin{subequations}
\begin{align}
&\MM^{(k)}_{\al_{i},\al_{j}}(x,n)=\MM^{(k)}_{\al_{j},\al_{i}}(-x,n)=
\left( \MM^{(k)}_{\al_{j},\al_{i}}(x,n) \right)^{*},
\qquad i,j=1,\ldots,d,
\label{eq:3.11a} \\
&\MM^{(k)}_{\be,\be}(x,n)=\MM^{(k)}_{\be,\be}(-x,n)=
\left( \MM^{(k)}_{\be,\be}(x,n) \right)^{*},
\label{eq:3.11b} \\
&\MM^{(k)}_{\al_{i},\be}(x,n)=-\MM^{(k)}_{\be,\al_{i}}(-x,n)=
- \left( \MM^{(k)}_{\be,\al_{i}}(x,n) \right)^{*},
\qquad i=1,\ldots,d ,
\label{eq:3.11c}
\end{align}
\label{eq:3.11}
\end{subequations}
\vskip-.5truecm \noindent  
where $*$ denotes complex conjugation.
\end{lemma}

\prova
Let us start from \eqref{eq:3.11a} --- in fact \eqref{eq:3.11b} can be obtained reasoning in the same way.
Given any $T\in\gotS^{k}_{n,\al_{i},\al_{j}}$ let $T'\in\gotS^{k}_{n,\al_{j},\al_{i}}$ be obtained
from $T$ by considering $\ell_{T},\ell'_{T}$ as entering and exiting lines, respectively, and reversing
the orientation of the lines in $\calP_{T}$. Denote by $N(\calP_{T})$ the set
of nodes in $N(T)$ connected by the lines in $\calP_{T}$. The node factors
of the nodes in $N(T)\setminus N(\calP_{T})$ and the
propagators of the lines outside $\calP_{T}$ do not change. Given $v\in
N(\calP_{T})$ let $\ell_{v},\ell'_{v}\in \calP_{T}\cup\{\ell_{T},\ell'_{T}\}$
be the lines exiting and entering $v$, respectively. If $h_{\ell_{v}}=
h_{\ell'_{v}}=\be$ or $h_{\ell_{v}},h_{\ell'_{v}}\in\{\al_{1},\ldots,\al_{d}\}$
then $\calF_{v}$ does not change when considering $v$ as a node in $T'$.
If $h_{\ell_{v}}=\be$ while $h_{\ell'_{v}}\in\{\al_{1},\ldots,\al_{d}\}$ or
vice versa, the node factor $\calF_{v}$ changes its sign when considering
$v$ as a node in $T'$. Since both $h_{\ell_{T}},
h_{\ell'_{T}}\in\{\al_{1},\ldots,\al_{d}\}$, then the number of nodes
in $N(\calP_{T})$ whose node factor changes sign must be even, so that
the overall product of such node factors does not change.
Finally if $\ell\in\calP_{T}$ one has $\nn_{\ell}=\nn_{\ell}^{0}+\nn_{\ell'_{T}}$
when considering it as a line in $L(T)$, while $\nn_{\ell}=-\nn_{\ell}^{0}+
\nn_{\ell'_{T'}}$ when considering it as a line in $\calP_{T'}$, so that,
computing at $\nn_{\ell'_{T'}}=-\nn_{\ell'_{T}}$, the propagators are equal
since they are even in their arguments. This proves the first equality in
\eqref{eq:3.11a}. Now let $T''\in\gotS^{k}_{n,\al_{j},\al_{i}}$ be obtained
from $T'$ by replacing the mode labels $\nn_{v}$ of the nodes in $N(T)$
with $-\nn_{v}$. The node factors are changed into their complex conjugated,
while (reasoning as before), when computing at $\nn_{\ell'_{T''}}=-\nn_{\ell'_{T}}$,
the propagators (which are real) do not change.

To prove \eqref{eq:3.11c} we reason as above, the only difference being that,
for $T\in\gotS^{k}_{n,\al_{j},\be}$, the numer of nodes in $N(\calP_{T})$
which change sing when considering them as nodes in $T'$ is odd, and
hence the overall product of the node factors change its sign.\EP

\begin{rmk}\label{rmk:3.9}
\emph{
 From Lemma \ref{lem:3.8} it follows that for all $k\ge1$ and all
$n\ge0$ one has
\begin{subequations}
\begin{align}
&\partial_{x}\MM^{(k)}_{\be,\be}(0,n)=0,
\nonumber \\ 
&\partial_{x}\MM^{(k)}_{\al_{i},\be}(0,n)=
- \left( \partial_{x}\MM^{(k)}_{\be,\al_{i}}(0,n) \right)^{*},
\qquad i=1,\ldots,d.
\nonumber 
\end{align}
\nonumber 
\end{subequations}
}
\end{rmk}

\begin{lemma}\label{lem:3.10}
For all $k\geq1$ one has
$\del_{x}\MM^{(k)}_{\al_{i},\al_{j}}(0,n)=0$ for $ i,j=1,\ldots,d$.
\end{lemma}

\prova
Given a cluster $T\in\gotS^{k}_{n,\al_{i},\al_{j}}$, with $i,j=1,\ldots,d$,
contributing to $M^{(k)}_{\al_{i},\al_{j}}(0,n)$ through (\ref{eq:3.9}), set
\begin{equation}\label{eq:3.13}
\del_{x} \Val_{T}(0):=\sum_{\ell\in\calP_{T}}\left(\prod_{v\in N(T)}\calF_{v}\right)
\left(\del_{x} \calG_{\ell}\prod_{\ell'\in L(T)\backslash\{\ell\}}\calG_{\ell'}\right),
\end{equation}
where the propagators have to be computed at $\oo\cdot\nn_{\ell'_{T}}=0$ and
\begin{equation} \nonumber 
\del_{x} \calG_{\ell}:=\frac{\Psi_{n_{\ell}}'(\oo\cdot\nn_{\ell}^{0})}{(\oo\cdot\nn_{\ell}^{0})^{2}}
-\frac{2\Psi_{n_{\ell}}(\oo\cdot\nn_{\ell}^{0})}{(\oo\cdot\nn_{\ell}^{0})^{3}},
\end{equation}
where $\Psi_{n}'$ denotes the derivative of $\Psi_{n}$ with respect to its argument.
Clearly $\partial_{x}\Val_{T}(0)$ is a contribution to
$\partial_{x}\MM^{(k)}_{\al_{i},\al_{j}}(0)$.

Now, the line $\ell$ divides $L(T)$ in two disjoint set of nodes 
$N_{1}$ and $N_{2}$ such that $\ell_{T}$ exits a node of $N_{1}$
and $\ell'_{T}$ enters a node in $N_{2}$. In other words if $\ell$ exits
a node $v$ one has $N_{2}=\{w\in N(T) : w\preceq v\}$
and $N_{1}=N(T)\setminus N_{2}$. Set
\begin{equation} \nonumber 
\nn_{1}=\sum_{v\in N_{1}}\nn_{v}; \qquad \nn_{2}=\sum_{v\in N_{2}}\nn_{v}.
\end{equation}
Since $T$ is a self-energy cluster one has $\nn_{1}+\nn_{2}=0$.
Then consider the family $\calF_{1}(T)$ of self-energy clusters obtained from 
$T$ by detaching the exiting line $\ell_{T}$ and reattaching it to all nodes
$w\in N_{1}$, and by detaching the entering line $\ell'_{T}$, then
reattaching it to all nodes $w\in N_{2}$. Consider also
a second family $\calF_{2}(T)$ of self-energy clusters obtained from $T$
by detaching the exiting line $\ell_{T}$ then reattaching it to all nodes 
$w\in N_{2}$ and by detaching the entering line $\ell'_{T}$ then reattaching
it to all nodes $w\in N_{1}$.

It can happen that, detaching $\ell_{T}$ from a node $w_{1}\in N_{1}$ and 
reattaching it to a node $w_{2}\in N(T)$, some node factors change their sign 
because some lines change their direction (see the proof of Lemma \ref{lem:3.8}).  
But, since $h_{\ell_{T}}=\al_{i}$ and $h_{\ell'_{T}}=\al_{j}$, the number of
changes of sign is even, so that 
the overall product of the node factors does not change its sign.
The shift of the lines $\ell_{T}$ and $\ell_{T}'$ also changes the combinatorial factors
of some node factors. However, if we group together all the self-energy clusters
in $\calF_{1}(T)$ with the two lines $\ell_{T}$ and $\ell_{T}'$ attached to the same nodes 
$v\in N_{1}$ and $w\in N_{2}$, respectively, we see that the corresponding values
differ from each other because of a factor $-\nn_{v}\nn_{w}$.
Reasoning in the same way we find that there are no changes 
of sign in the product of the node factors also in the construction
of the family $\calF_{2}(T)$. Moreover, for those lines that
change their direction after such shift operation, the momentum $\nn_{\ell}$ 
is replaced by $-\nn_{\ell}$ but no changes are produced 
in the propagators since they are even, except for the differentiated
propagator which can change sign: the sign changes for the self-energy clusters
in $\calF_{1}(T)$, while it remains the same for those in $\calF_{2}(T)$. 
Then by summing over all possible cluster in $\calF_{1}(T)$ 
we obtain $-\nn_{1}\nn_{2}$ times a common factor, while summing over all 
possible cluster in $\calF_{2}(T)$ we obtain $\nn_{1}\nn_{2}$ times the
same common factor, so that the overall sum gives zero.\EP

\begin{lemma}\label{lem:3.11}
For all $k\ge1$ one has
\begin{subequations}
\begin{align}
&[-\partial_{\aaa}f(\aaa,\be)]^{(k)}_{\vzero}=\vzero
\label{eq:3.14a} \\
&\MM^{(k)}_{\al_{i},h}(0)=0,
\qquad i=1,\ldots,d,\quad h=\al_{1},\ldots,\al_{d},\be
\label{eq:3.14b} 
\end{align}
\label{eq:3.14}
\end{subequations}
\vskip-.6truecm \noindent
\end{lemma}

\prova
We first prove \eqref{eq:3.14a}. Given $\theta \in
\Theta_{k,\vzero,\al_{j}}$, denote by $\gotF(\theta)$ the set of all possible
$\theta'\in \Theta_{k,\vzero,\al_{j}}$ which can be obtained from $\theta$
by detaching the root line $\ell_{\theta}$ and reattaching it to each node
$v\in N(\theta)$. The values of such trees differ from each other
because of a factor $\ii\nn_{v}$, where $v$ is the node which the
root line is attached to (again, as in the proof of Lemmas \ref{lem:3.8}
and \ref{lem:3.10}, there is an even number of nodes whose node factor
changes sign, and hence the overall product does not changes). But then, since
$\sum_{v\in N(\theta)}\nn_{v}=\vzero$, the sum over all such contributions is zero.
Moreover this holds identically in $\aaa_{0},\be_{0}$, therefore by Lemma \ref{lem:3.7} also
\eqref{eq:3.14b} follows.\EP

\begin{rmk}\label{rmk:3.12} 
\emph{
Identity \eqref{eq:3.14a} is formally equal to \eqref{eq:2.2c}: therefore
we proved that \eqref{eq:2.2c} formally holds. So, besides the
convergence of the series, we are left with \eqref{eq:2.2d} to be solved.
}
\end{rmk}

We can summarise the results above as follows.
Let us write
\begin{equation}\label{eq:3.15}
\MM^{(k)}_{u,e}(x,n)=\matL_{u,e}^{(k)}+ x \DD_{u,e}^{(k)}+
x^{2}\matD_{u,e}^{(k)}(x)+ \matR_{u,e}^{(k)}(x,n),
\end{equation}
with
\begin{equation}\label{eq:3.16}
\begin{aligned}
& \matL_{u,e}^{(k)}:=\MM^{(k)}_{u,e}(0), \qquad
\DD_{u,e}^{(k)}:=\partial_{x}\MM^{(k)}_{u,e}(0),\qquad
\matD_{u,e}^{k}(x,n):=\int_{0}^{1}{\rm d}\tau\,(1-\tau)\partial^{2}_{x}
\MM^{(k)}_{u,e}(\tau x),\\
& \matR_{u,e}^{(k)}(x,n):=\MM^{(k)}_{u,e}(x,n) - \MM^{(k)}_{u,e}(x).
\end{aligned}
\end{equation}
Then we have
\begin{equation}\label{eq:3.17}
\begin{aligned}
&\MM^{(k)}_{\al_{i},\al_{j}}(x,n)=
x^{2}\matD_{\al_{i},\al_{j}}^{(k)}(x)+ \matR_{\al_{i},\al_{j}}^{(k)}(x,n),
\qquad i,j=1,\ldots,d, \\
&\MM^{(k)}_{\al_{i},\be}(x,n)= x \DD_{\al_{i},\be}^{(k)}+
x^{2}\matD_{\al_{i},\be}^{(k)}(x)+ \matR_{\al_{i},\be}^{(k)}(x,n),
\qquad i=1,\ldots,d, \\
&\MM^{(k)}_{\be,\al_{i}}(x,n)= - x \big( \DD_{\al_{i},\be}^{(k)} \big)^{*}+
x^{2}\matD_{\be,\al_{i}}^{(k)}(x)+ \matR_{\be,\al_{i}}^{(k)}(x,n),
\qquad i=1,\ldots,d, \\
&\MM^{(k)}_{\be,\be}(x,n)=\matL_{\be,\be}^{(k)}+
x^{2}\matD_{\be,\be}^{(k)}(x)+ \matR_{\be,\be}^{(k)}(x,n).
\end{aligned}
\end{equation}
In other words, if we could ingore the ``rest'' $\matR_{u,e}^{(k)}(x,n)$, we would
obtain a gain factor proportional to $x^{2}$ for the self-energies with
$u,e\in\{\al_{1},\ldots,\al_{d}\}$, a gain proportional to $x$ for $u=\al_{1},
\ldots,\al_{d}$ and $e=\be$ (or vice versa) and no gain for $u=e=\be$ (but in
the latter case no factor proportional to $x$ would appear). This suggests us
that if $\matL_{\be,\be}^{(k)}\equiv 0$
and $\DD_{\al_{i},\be}^{(k)}\equiv0$ for all $i=1,\ldots,d$ and all $k\ge1$,
we would obtain a gain proportional to $x^{2}$ for any self-energy (provided
the ``rest'' is small) and this should imply the convergence of the power series.

\begin{cond}\label{cond:1}
For all $k\ge1$ one has $\matL_{\be,\be}^{(k)}\equiv 0$
and $\DD_{\al_{i},\be}^{(k)}\equiv0$ for all $i=1,\ldots,d$.
\end{cond}

\begin{lemma}\label{lem:3.13}
Assume Condition \ref{cond:1}.
Then for all $h,h'\in\{\al_{1},\ldots,\al_{d},\be\}$
and for any $(\aaa_{0},\be_{0})\in\TTT^{d+1}$ one has
$|\MM^{(k)}_{h,h'}(x,n)|\le C^{k}x^{2}$,
for some positive constant $C$.
\end{lemma}

The proof of the result above essentially follows the lines of the proof
of Lemma 6.6 in \cite{CG2}.

\begin{rmk}\label{rmk:3.14}
\emph{
One can prove also that, setting $\DD_{\aaa,\be}^{(k)}=(\DD_{\al_{1},\be}^{(k)},\ldots,
\DD_{\al_{d},\be}^{(k)})$, one has
\begin{equation}\nonumber
\oo\cdot\DD_{\aaa,\be}^{(k)}=2\ii(k-1)[\partial_{\be}f(\aaa,\be)]^{(k)}_{\vzero},
\end{equation}
for all $k\ge1$. We shall not give the proof of the identity above
since it will not be used here.
}
\end{rmk}

If Condition \ref{cond:1} is satisfied,
Lemma \ref{lem:3.13} implies the convergence of the series
(\ref{eq:2.4}) for $\e$ small enough:
the argument is the same as after Lemma 6.6 in \cite{CG2}.
Moreover, by Lemma \ref{lem:3.7}, the assumption
$\matL_{\be,\be}^{(k)}\equiv 0$ reads
$[\partial_{\be}f(\aaa,\be)]^{(k)}_{\vzero}\equiv \hbox{const}$. Due
to the variational nature of the Hamilton equation, $[\partial_{\be}f(\aaa,\be)
]^{(k)}_{\vzero}$ is the $\be_{0}$-derivative of the $k$-th order of the time
average of the Lagrangian $\g^{(k)}$ (which is analytic and periodic) computed
along a solution of the range equation (one can reason as in \cite{CG3}).
This implies $[\partial_{\be}f(\aaa,\be)]^{(k)}_{\vzero}\equiv0$, so that also
(\ref{eq:2.2d}) holds for any $\be_{0}\in\TTT^{d}$. Therefore, at least in
the particular case that Condition \ref{cond:1} holds, we provided
a quasi-periodic solution to the equation (\ref{eq:1.2}) as a convergent
power series in $\e$. Note that in such a case the initial phase $\be_{0}$ remains
arbitrary, so that the full $(d+1)$-resonant unperturbed torus persists.

\zerarcounters 
\section{Resummation of the formal expansion} 
\label{sec:4} 

In Section \ref{sec:3} we have seen how to deal with the ``completely
degenerate case'' of Condition \ref{cond:1}, which yields infinitely many identities.
If these identities do not hold we are not able to prove the convergence
of the series \eqref{eq:2.4}. Now we shall see how to deal with such a case.

\subsection{Renormalised trees}\label{sec:4.1}

As seen in Section \ref{sec:3.2} all the obstruction to the convergence
are due to the presence of self-energy clusters. Now we shall perform a
different tree expansion with respect to the one performed in Section
\ref{sec:3.1} in order to deal with this problem.

More precisely, we modify the tree expansion envisaged in Section \ref{sec:3.1}
as follows. Given a tree $\theta$ we associate with each node $v\in
N(\theta)$ a mode label and a component label as in Section \ref{sec:3.1};
with each line $\ell\in L(\theta)$ we associate a momentum label
as in Section \ref{sec:3.1} and a \emph{pair} of component labels
$(e_{\ell},u_{\ell})\in\{\al_{1},\ldots,\al_{d}\}$ with the constraint that
$u_{\ell_{v}}=h_{v}$. We shall call $e_{\ell}$ and $u_{\ell}$ the
$e$-\emph{component} and the $u$-\emph{component} of $\ell$, respectively.
We denote by $p_{v}$ and $q_{v}$ the number of lines with $e$-component
$\al_{j}$ for some $j=1,\ldots,d$ and $\be$ entering $v$, respectively,
and set $s_{v}=
p_{v}+q_{v}$. We still impose the conservation law \eqref{eq:3.1}.
We do not change the definition of cluster, while from now on a self-energy
cluster is a cluster $T$ with only one entering line $\ell'_{T}$ and one
exiting line $\ell_{T}$ such that $\nn_{\ell_{T}}=\nn_{\ell'_{T}}$, i.e.
we drop the constraint (ii) from the definition of self-energy cluster
given in Section \ref{sec:3.1}.

A \emph{renormalised tree} is a tree in which no self-energy cluster appears.
Analogously a \emph{renormalised subgraph} of a tree is a subgraph
$S$ of a tree $\theta$ such that $S$ does not contains any self-energy cluster.

Given a renormalised tree we call \emph{total momentum} and
\emph{total component} the momentum and the $e$-component associated with
the root line. We denote by $\Theta^{\RR}_{k,\nn,h}$ the set of all
renormalised trees with order $k$ total momentum $\nn$ and total
component $h$, and by $\gotR_{n,u,e}$ the set of renormalised self-energy
clusters on scale $n$ such that $u_{\ell_{T}}=u$ and $e_{\ell'_{T}}=e$.

Given $\theta\in\Theta^{\RR}_{k,\nn,h}$ we associate with each $v\in N(\theta)$
a node factor $\calF_{v}$ defined as in \eqref{eq:3.2} and with each
$\ell\in L(\theta)$ a propagator $\calG_{\ell}$ defined as follows.
First of all, given a $(d+1)\times(d+1)$ matrix $A$ with entries $A_{h,h'}$,
for $h,h'\in\{\al_{1},\ldots,\al_{d},\be\}$, we denote by
$A_{\aaa,\aaa}$ the $d\times d$ matrix with entries
$(A_{\aaa,\aaa})_{i,j}:=A_{\al_{i},\al_{j}}$, for $i,j=1,\ldots, d$,
by $A_{\aaa,\be}$ the vector with components
$(A_{\aaa,\be})_{i}:=A_{\al_{i},\be}$, for $i=1,\ldots, d$, and
by $A_{\be,\aaa} $ the vector with components
$(A_{\be,\aaa})_{j}:=A_{\be,\al_{j}}(x)$, for $j=1,\ldots, d$;
with a slight abuse of notation we denote in the same way both column and row vectors.
Then we define recursively the propagator of the line $\ell$ as
$\calG_{\ell}:=\calG_{e_{\ell},u_{\ell}}^{[n_{\ell}]}(\oo\cdot\nn_{\ell})$, with
\begin{equation}\label{eq:4.1}
\calG^{[n]} (x)=
\begin{pmatrix}
\calG^{[n]}_{\aaa,\aaa}(x) & \calG^{[n]}_{\aaa,\be}(x) \cr
\calG^{[n]}_{\be,\aaa}(x) & \calG^{[n]}_{\be,\be}(x)
\end{pmatrix}:=
\Psi_{n}(x)\left(x^{2}\uno - \MM^{[n-1]}(x)\right)^{-1} ,
\end{equation}
where $\Psi_{n}$ is defined as in Section \ref{sec:3.1},
$\uno$ is the $(d+1)\times(d+1)$ identity matrix and
\begin{equation}\label{eq:4.2}
\MM^{[n-1]}(x):=\sum_{q=-1}^{n-1}\chi_{q}(x)M^{[q]}(x),
\end{equation}
with $\chi_{q}$ defined as in Section \ref{sec:3.1} and
\begin{equation}\label{eq:4.3}
M^{[n]}(x)=\left(\begin{matrix}
M_{\aaa,\aaa}^{[n]}(x) & M_{\aaa,\be}^{[n]}(x) \\
M_{\be,\aaa}^{[n]}(x) & M_{\be,\be}^{[n]}(x)\end{matrix}\right) ,
\;\mbox{ with }\;
M_{u,e}^{[n]}(x):=\sum_{T\in\gotR_{n,u,e}}\e^{k(T)}\Val_{T}(x),
\quad n\ge-1,
\end{equation}
and
\begin{equation}\label{eq:4.4}
\Val_{T}(x):=\left(\prod_{v\in N(T)}\calF_{v}\right)
\left(\prod_{\ell\in L(T)}
\calG_{\ell}\right)
\end{equation}
is the {\emph{renormalised value}} of $T$.

Set $\MM:=\{\MM^{[n]}(x)\}_{n\ge-1}$. We call \emph{self-energies} the
matrices $\MM^{[n]}(x)$.

\begin{rmk}\label{rmk:4.1}
\emph{
By construction $\calG^{[n]}(x)$ depends also on $\e$ and $\be_{0}$,
even though we are not making explicit such a dependence;
it does not depend on $\aaa_{0}$ because 
$f_{\nn}(\aaa_{0},\be_{0})={\rm e}^{{\rm i}\nn\cdot\aaa_{0}} \hat f_{\nn}(\be_{0})$
and $\sum_{v\in N(T)}\nn_{v}=0$ for any self-energy cluster $T$.
The last comment applies also to the quantities
$\FF^{[k]}(\e,\be_0)$ and $G^{[k]}(\e,\be_0)$ introduced in (\ref{eq:4.5}) below.
}
\end{rmk}

Setting also $\calG^{[-1]}=\uno$, for any renormalised subgraph $S$
of any $\theta\in\Theta^{\RR}_{k,\nn,h}$ we define the \emph{renormalised value} of $S$ as
in \eqref{eq:3.5}, but with the new definition for the propagators.

Set $\Theta^{\RR}_{k,\nn,\aaa}:=\Theta^{\RR}_{k,\nn,\al_1}\times\ldots\times
\Theta^{\RR}_{k,\nn,\al_d}$ and for any
$\ttheta=(\theta_1,\ldots,\theta_d)\in\Theta^{\RR}_{k,\nn,\aaa}$
denote $\VVal(\ttheta):=\left(\Val(\theta_1),\ldots,\Val(\theta_d)\right)$.
Then define (formally)
\begin{equation}\label{eq:4.5}\begin{aligned}
&\aaaa_{\nn}^{[k]}(\e,\aaa_0,\be_0):=\!\!\!\!\!\!
\sum_{\ttheta\in\Theta_{k,\nn,\aaa}^{\RR}} \!\!\! \VVal(\ttheta), \qquad
b_{\nn}^{[k]}(\e,\aaa_{0},\be_{0}):=\!\!\!\!\!\!
\sum_{\theta\in\Theta_{k,\nn,\be}^{\RR}} \!\!\! \Val(\theta), \qquad
\nn\neq\vzero, \\
&\FF^{[k]}(\e,\be_0):=\!\!\!\!\!\!
\sum_{\ttheta\in\Theta_{k+1,\vzero,\aaa}^{\RR}} \!\!\! \VVal(\ttheta), \qquad
G^{[k]}(\e,\be_0):=\!\!\!\!\!\!
\sum_{\theta\in\Theta_{k+1,\vzero,\be}^{\RR}} \!\!\! \Val(\theta).
\end{aligned}
\end{equation}
Finally set (again formally)
\begin{subequations}
\begin{align}
&\aaaa^{\RR}(t;\e,\aaa_0,\be_0):=\sum_{k\geq1}\e^{k}\sum_{\nn\in\ZZZ_{*}^{d}}
{\rm e}^{\ii\nn\cdot\oo t}\aaaa_{\nn}^{[k]}(\e,\aaa_0,\be_0),
\label{eq:4.6a}\\
&b^{\RR}(t;\e,\aaa_0,\be_0):=\sum_{k\geq1}\e^{k}\sum_{\nn\in\ZZZ_{*}^{d}}
{\rm e}^{\ii\nn\cdot\oo t}b_{\nn}^{[k]}(\e,\aaa_0,\be_0),
\label{eq:4.6b}\\
&\FF^{\RR}(\e,\be_0):=\sum_{k\geq0}\e^{k}\FF^{[k]}(\e,\be_0),
\label{eq:4.6c}\\
&G^{\RR}(\e,\be_0):=\sum_{k\geq0}\e^{k}G^{[k]}(\e,\be_0) ,
\label{eq:4.6d}
\end{align}
\label{eq:4.6}
\end{subequations}
\vskip-.1truecm \noindent 
and define
\begin{equation}\label{eq:4.7}
\begin{aligned}
&\aaa^{\RR}(t;\e,\aaa_0,\be_0)=\aaa_0+\oo t+\aaaa^{\RR}(t;\e,\aaa_0,\be_0),\\
&\be^{\RR}(t;\e,\aaa_0,\be_0)=\be_0+b^{\RR}(t;\e,\aaa_0,\be_0).
\end{aligned}
\end{equation}
The series (\ref{eq:4.6}) will be called \emph{resummed series}, the term
``resummed'' coming from the fact that if we formally expand
\eqref{eq:4.6} in powers of $\e$ then we get \eqref{eq:2.4}, as is easy to check.

For any renormalised subgraph $S$ of any tree $\theta$ we denote by
$\gotN_{n}(S)$ the number of lines on scale $\ge n$ in $S$ and define
$K(S)$ as in \eqref{eq:3.7}.
Then we have the following results which are the counterparts of Lemmas
\ref{lem:3.3} and \ref{lem:3.4}, respectively, for renormalised trees.

\begin{lemma}\label{lem:4.2}
For any $h\in\{\al_{1},\ldots,\al_{d},\be\}$, $\nn\in\ZZZ^{d}$, $k\ge 1$
and for any $\theta\in \Theta_{k,\nn,h}^{\RR}$ such
that $\Val(\theta)\ne 0$, one has $\gotN_{n}(\theta)\le 
2^{-(m_{n}-2)}K(\theta)$ for all $n\ge0$.
\end{lemma}

\begin{lemma}\label{lem:4.3}
For any $e,u\in\{\al_{1},\ldots,\al_{d},\be\}$, $n \ge0$ and for any
$T\in\gotR_{n,u,e}$ such that
$\Val_{T}(x)\ne 0$, one has $K(T) > 2^{m_{n}-1}$ and
$\gotN_{p}(T)\le 2^{-(m_{p}-2)} K(T)$ for $0\le p\le n$.
\end{lemma}

The two results above can be proved as Lemmas
4.1 and 4.2 in \cite{CG2}, respectively.

\subsection{A suitable assumption: bounds}\label{sec:4.2}

Here we shall see that, under the assumption that the propagators
$\calG^{[n]}_{e,u}(\oo\cdot\nn)$ are
bounded proportionally to $1/|\oo\cdot\nn|^{c}$ for some constant $c$, the series
(\ref{eq:4.6}) converge and solve the range equations
(\ref{eq:2.2a}) and \eqref{eq:2.2b}: the key point is that now self-energy
clusters (and hence resonant lines) are not allowed and hence a result of
that kind is expected. Then,
in what follows, we shall see that the assumption is justified
at least along a curve $\be_{0}(\e)$ where
also the bifurcation equations (\ref{eq:2.2c}) and (\ref{eq:2.2d}) are satisfied.

Define $\|\cdot\|$ as an algebraic matrix norm (i.e. a norm which verifies
$\| AB\|\leq\|A\| \|B\|$ for all matrices $A$ and $B$);
for instance $\|\cdot\|$ can be the uniform norm.

\begin{defi}\label{def:4.4} 
We shall say that $\MM$ satisfies Property 1 if
there are positive constants $c_{1}$ and $c_{2}$ such that
$$
\left\| \calG^{[n]}(x)\right\|\leq\frac{c_{1}}{|x|^{c_{2}}}
$$
for all $n\ge 0$. Call $\SSSS:=\left\{(\e,\be_{0})\in\RRR\times\TTT : 
\hbox{Property 1 holds} \right\}$.
\end{defi}

\begin{defi}\label{def:4.5}
We shall say that $\MM$ satisfies Property 1-$p$ if
there are positive constants $c_{1}$ and $c_{2}$ such that
$$
\left\| \calG^{[n]}(x)\right\|\leq\frac{c_{1}}{|x|^{c_{2}}},
$$
for $0\le n\le p$. Call $\SSSS_{p}:=\left\{(\e,\be_{0})\in\RRR\times\TTT : 
\hbox{Property 1-$p$ holds}\right\}$.
\end{defi}

\begin{lemma}\label{lem:4.6} 
Assume $(\e,\be_{0})\in\SSSS_{p}$. Then, for
$0\le n\le p$ and $\e$ small enough, the self-energies are well defined
and one has 
\begin{equation} \nonumber 
\left| \partial_{x}^{j}M^{[n]}_{u,e}(x) \right|\le  \e^{2} \, K_{j}
{\rm e}^{-\ol{K}_{j}2^{m_{n}}},\qquad j=0,1,2,
\end{equation}
for some positive constants $K_{0},\ol{K}_{0},K_{1},\ol{K}_{1},K_{2}$ and $\ol{K}_{2}$.
\end{lemma}

The proof is essentially the same as the proof of Lemma
4.8 in \cite{CG2} and Lemma 4.3 in \cite{CG3}.
In particular we need Remark \ref{rmk:3.6} 
when bounding the derivatives.

\begin{rmk}\label{rmk:4.7}
\emph{
If $\MM$ satisfies Property 1-$p$
the matrices $\MM^{[n]}(x)$ and $\calG^{[n]}(x)$ are well defined
for all $-1\le n\le p$. In particular there exists $\g_{0}>0$ such that
$|\calG^{[n]}_{e,u}(x)|\le \g_{0}\,\al_{m_{n}}(\oo)^{-c_{2}}$ for all $0\le n\le p$.
If $\MM$ satisfies Property $1$, the same considerations
apply for all $n\ge 0$.
}
\end{rmk}

\begin{lemma}\label{lem:4.8}
Assume $(\e,\be_{0})\in\SSSS_{p}$. Then, for
$0\le n\le p$ and $\e$ small enough, one has 
\begin{equation} \nonumber 
\left| M^{[n]}_{u,e}(x) - M^{[n]}_{u,e}(0) - x \,
\partial_{x} M^{[n]}_{u,e}(0) \right|
\le  \e^{2} \,  {K}_{3} {\rm e}^{-\ol{K}_{3} 2^{m_{n}}} x^{2}
\end{equation}
for some positive constants $K_{3}$ and $\ol{K}_{3}$.
\end{lemma}

The proof is essentially the same as the proof of Lemma 4.6 in \cite{CG3}.

\begin{lemma}\label{lem:4.9} 
Assume $(\e,\be_{0})\in\SSSS$. Then 
the series (\ref{eq:4.6}), with the coefficients 
given by (\ref{eq:4.5}), converge for $\e$ small enough. 
\end{lemma} 

The proof is essentially the same as the proof of Lemma 4.5 in \cite{CG2}
and Lemma 4.9 in \cite{CG3}.

\begin{lemma}\label{lem:4.10} 
Assume $(\e,\be_{0})\in\SSSS$. Then for $\e$ small enough the
function (\ref{eq:4.6a}) and (\ref{eq:4.6b}) solve the
range equations (\ref{eq:2.2a}) and (\ref{eq:2.2b}), respectively. 
\end{lemma} 

The proof is essentially the same as the proof of Lemma 4.6 in \cite{CG2}
and Lemma 4.10 in \cite{CG3}.

\subsection{A suitable assumption: symmetries}\label{sec:4.3}

Here we shall prove that, under the assumptions that $\MM$ satisfies
Property 1-$p$, there are suitable symmetries for the self-energy
clusters: such symmetries are the counterpart of those founded
in Section \ref{sec:3.3} for the formal expansion. Property 1-$p$ is
assumed only because, under such assumption, all the quantities are well defined.

\begin{lemma}\label{lem:4.11}
Let $\gotB_{n}$ the set of $B:\RRR\to {\rm GL}(n,\CCC)$ such that
\begin{subequations} \nonumber
\begin{align}
& B_{i,j}(-x)=B_{j,i}(x) , \qquad i,j=1,\ldots,n-1, 
\qquad B_{n,n}(-x)=B_{n,n}(x) \\
& B_{n,i}(-x)=-B_{i,n}(x) , \qquad i=1,\ldots,n-1 .
\end{align}
\end{subequations}
\vskip-.1truecm \noindent
Then if $B \in \gotB_{n}$ also $B^{-1}\in\gotB_{n}$.
\end{lemma}

\prova
If $B\in\gotB_{n}$ define the matrix $A$ by setting
\begin{subequations} \nonumber
\begin{align}
& B_{i,j}(x)=A_{i,j}(x) , \qquad i,j=1,\ldots,n-1, \qquad
B_{n,n}(x)=A_{n,n}(x) , \\
& B_{n,i}(x)=x \, A_{n,i}(x) \hbox{ and } B_{i,n}(x)=x \, A_{i,n}(x) , \qquad i=1,\ldots,n-1, 
\end{align}
\end{subequations}
\vskip-.1truecm \noindent
so that $A^{T}(-x)=A(x)$. Denote also by $C_{i,j}(x)$ the cofactor of the
entry $A_{i,j}(x)$ for $i,j=1,\ldots,n$. By construction
$C_{i,j}(-x)=C_{j,i}(x)$ for $i,j=1,\ldots,n$: then
\begin{equation}\nonumber
 \begin{aligned}
\det B(x) &=  (-1)^{n-1} x^{2} \left[ A_{n,1}(x) \, C_{n,1}(x) - \ldots +
(-1)^{n-2} A_{n,n-1}(x)\,C_{n,n-1} (x) \right]
+ A_{n,n}(x) \, C_{n,n}(x)   \\
&= x^{2} \det A(x) + \left( 1 - x^{2} \right) A_{n,n}(x) \, C_{n,n}(x) ,
\end{aligned}
\end{equation}
so that $\det B(-x)=\det B(x)$. By noting that
\begin{equation} \nonumber
\left( B^{-1}(x)\right)_{j,i} = \frac{1}{\det B(x)} \begin{cases}
x^{2} C_{i,j}(x)+(1-x^{2})D_{i,j}(x) , & i,j=1,\ldots,n-1 , \\
x \, C_{i,j}(x) , & i=n \hbox{ and } j=1,\ldots,n-1 , \\
x \, C_{i,j}(x) , & i=1,\ldots,n-1  \hbox{ and }  j=n , \\
C_{i,j}(x) , & i,j=n , \end{cases}
\end{equation}
where $D_{i,j}(x)$ is the cofactor of $A_{i,j}(x)$ seen as entry of the $(n-1)\times(n-1)$ matrix
obtained from $A(x)$ by deleting its $n$-th row and $n$-th column,
the assertion follows.\EP

\begin{lemma}\label{lem:4.12}
Assume $(\e,\be_{0})\in\SSSS_{p}$. Then for all $-1\le n\leq p$ one has
\begin{subequations} 
\begin{align}
\left( \MM^{[n]}_{\aaa,\aaa}(x)\right)^{T}
& =\MM^{[n]}_{\aaa,\aaa}(-x)
=\left(\MM^{[n]}_{\aaa,\aaa}(x)\right)^{*},
\label{eq:4.10a} \\ 
\MM^{[n]}_{\be,\be}(x)
& =\MM^{[n]}_{\be,\be}(-x)
=\left(\MM^{[n]}_{\be,\be}(x)\right)^{*},
\label{eq:4.10b} \\
\MM^{[n]}_{\aaa,\be}(x)
& = -\MM^{[n]}_{\be,\aaa}(-x)=
-\left(\MM^{[n]}_{\be,\aaa}(x)\right)^{*},
\label{eq:4.10c}
\end{align}
\label{eq:4.10}
\end{subequations}
\vskip-.4truecm \noindent
where $*$ denotes complex conjugation.
\end{lemma}

\prova
We shall proceed by induction on $n$. First of all note that for $n=-1$
(\ref{eq:4.10}) trivially holds, since
\begin{equation} \label{eq:04.16}
\MM^{[-1]}(x) = 
\begin{pmatrix}
0_{d} & \vzero \cr \vzero & \e\partial_{\be}^{2} f_{\vzero}
\end{pmatrix},
\end{equation}
where $0_{d}$ is the $d\times d$ null matrix. Assume than that
\eqref{eq:4.10} hold for all $-1\le n'<n$ and let us start from
the first equality in \eqref{eq:4.10a}. Given any $T\in\gotR_{n,\al_{i},\al_{j}}$
let $T'$ be obtained from $T$ by reversing the orientation of the
lines along $\calP_{T}\cup \{\ell_{T},\ell'_{T}\}$.
Denote by $N(\calP_{T})$ the set
of nodes in $N(T)$ connected by the lines in $\calP_{T}$.
The node factors of the nodes in $N(T)\setminus N(\calP_{T})$ and the
propagators on the lines outside $\calP_{T}$ do not change when
considering them as nodes and lines in $T'$. Given $v\in
N(\calP_{T})$ let $\ell_{v},\ell'_{v}\in \calP_{T}\cup\{\ell_{T},\ell'_{T}\}$
be the lines exiting and entering $v$, respectively. If $u_{\ell_{v}}=
e_{\ell'_{v}}=\be$ or $u_{\ell_{v}},e_{\ell'_{v}}\in\{\al_{1},\ldots,\al_{d}\}$
then $\calF_{v}$ does not change when considering $v$ as a node in $T'$.
If $u_{\ell_{v}}=\be$ while $e_{\ell'_{v}}\in\{\al_{1},\ldots,\al_{d}\}$ or
vice versa, the node factor $\calF_{v}$ changes its sign when considering
$v$ as a node in $T'$. Now, given $\ell\in\calP_{T}$ we compute the propagator
associated with $\ell$ at $x_{\ell}:=\oo\cdot\nn_{\ell}=\oo\cdot\nn_{\ell}^{0}+x$
and we obtain $\calG_{\ell}=\Psi_{n_{\ell}}(x_{\ell})
(x_{\ell}^{2}\uno-\MM^{[n_{\ell}-1]}(x_{\ell}))_{e_{\ell},u_{\ell}}^{-1}$; when considering
$\ell$ as a line in $T'$, if we set $\oo\cdot\nn_{\ell'_{T'}}=-x$, then
the momentum of $\ell$ changes sign and hence the propagator becomes
$\Psi_{n_{\ell}}(-x_{\ell})
((-x_{\ell})^{2}\uno-\MM^{[n_{\ell}-1]}(-x_{\ell}))_{u_{\ell},e_{\ell}}^{-1}$: thanks to the
inductive hypothesis and Lemma \ref{lem:4.11}, if $e_{\ell}=u_{\ell}=
\be$ or $e_{\ell}=u_{\ell}\in\{\al_{1},\ldots,\al_{d}\}$ the propagator
does not change when considering $\ell$ as a line in $T'$, otherwise it
changes its sign. Let $h_{0},\ldots,h_{2|\calP_{T}|+1}$
be such that $h_{0}=u_{\ell_{T}}$, $\{h_{1},\ldots,h_{2|\calP_{T}|}\}$ is the ordered set of
the components of the lines in $\calP_{T}$ and $h_{2|\calP_{T}|+1}=e_{\ell'_{T}}$.
Note that there is a change of sign (in the node factor or in the propagator)
corresponding to each ordered pair $h_{r},h_{r+1}$ such that either
$h_{r}=\al_{i}$ for some $i=1,\ldots,d$ and $h_{r+1}=\be$ or vice versa.
Since $h_{0}=\al_{i}$ and $h_{2|\calP_{T}|+1}=\al_{j}$ the number of changes of
sign is even and therefore the overall product does not change.
This proves the first equality in \eqref{eq:4.10a}.
The first equality in \eqref{eq:4.10b} can be proved in the same way.

Now let $T''$ be the self-energy cluster obtained from $T'$ by replacing
the mode labels $\nn_{v}$ of the nodes in $N(T')$ with $-\nn_{v}$.
The node factors are changed into their complex conjugated and, thanks
to the inductive hypothesis, when computing at $\nn_{\ell'_{T''}}=
-\nn_{\ell'_{T'}}$, also the propagators are changed into their complex
conjugated. Hence also the second equality in \eqref{eq:4.10a} is proved.
Again analogous considerations lead to the second equality in (\ref{eq:4.10b}).

To prove \eqref{eq:4.10c} one can reason
in the same way, the only difference being that for $T\in
\gotR_{n,\al_{i},\be}$ the number of changes of sign of the propagators of
the lines in $\calP_{T}$ and of the node factors of the nodes in $N(\calP_{T})$ is
odd. This implies the change of sing in the first equality in \eqref{eq:4.10c}.\EP

\begin{lemma}\label{lem:4.13}
Assume $(\e,\be_{0})\in\SSSS_{p}$. Then one has for $-1\le n\le p$
\begin{subequations}
\begin{align}
&\MM^{[n]}_{\aaa,\aaa}(x)=O(\e^{2}x^2),
\label{eq:4.12a}\\
&\MM^{[n]}_{\be,\aaa}(x)=O(\e^{2}x),
\label{eq:4.12b}\\
&\MM^{[n]}_{\aaa,\be}(x)=O(\e^{2}x),
\label{eq:4.12c}\\
&\MM^{[n]}_{\be,\be}(x)=\MM^{[n]}_{\be,\be}(0)+O(\e^{2}x^{2}),
\label{eq:4.12d}
\end{align}
\label{eq:4.12}
\end{subequations}
\vskip-.4truecm \noindent
where $\MM^{[n]}_{\be,\be}(0)=O(\e)$.
\end{lemma}

\prova
Let us start from the proof of (\ref{eq:4.12a}). First of all we shall
show that  $\sum_{T\in\gotR_{n,u,e}}\Val_{T}(0)=0$ where $(u,e)\in
\{\al_1,\ldots,\al_d\}^{2}$.
Given a self-energy cluster $T\in\gotR_{n,\al_i,\al_j}$ for $i,j=1,\ldots,d$ 
consider all the self-energy cluster which can be obtained from $T$ by
detaching the entering line $\ell'_{T}$ and reattaching it to each node $v\in N(T)$. 
After such operation $\Val_{T}(0)$ changes by a factor $(\ii\nn_{v})$ 
if $v$ is the node which the entering line is attached to, 
while the other node factors and propagators do not change
(the combinatorial factors can be discussed as along the proof of
Lemma \ref{lem:3.10}).  The sum of all clusters values is zero because
$\sum_{v\in N(T)}\nn_{v}=0$. This implies
$\MM^{[n]}_{\aaa,\aaa}(0)=0_{d}$
(see the beginning of the proof of Lemma \ref{lem:4.12} for notation).

Now let us write $\partial_{x}\Val_{T}(0)$ as in \eqref{eq:3.13},
where again the propagators have to be computed at $\oo\cdot\nn_{\ell_{T}'}=0$, but now
\begin{equation} \nonumber 
\partial_{x}\calG_{\ell}:=\left.\frac{{\rm d}}{{\rm d} x}
\calG^{[n_{\ell}]}_{e_{\ell},u_{\ell}}(\oo\cdot\nn_{\ell}^{0}+x)\right|_{x=0}.
\end{equation}
The line $\ell$ divides $L(T)$ in two disjoint set of nodes $N_1$ and $N_2$ 
such that $\ell_{T}$ exits a node in $N_{1}$ and $\ell'_{T}$ enters a node
in $N_{2}$.  In other words $N_{2}=\{w\in N(T) : w\prec \ell\}$ and
$N_{1}=N(T)\setminus N_{2}$. Set
 \begin{equation}\label{eq:4.13}
 \nn_{1}=\sum_{v\in N_{1}}\nn_{v}; 
 \qquad \nn_{2}=\sum_{v\in N_{2}}\nn_{v}.
 \end{equation}
Since $T$ is a self-energy cluster one has $\nn_{1}+\nn_{2}=0$. 
Now consider the family $\calF_{1}(T)$ of self-energy cluster obtained
from $T$  by detaching the exiting line $\ell_{T}$ and reattaching it
to all nodes $w\in N_{1}$, and by detaching the entering line $\ell'_{T}$
and reattaching it to all nodes $w\in N_{2}$. Consider also the family
$\calF_{2}(T)$ obtained from $T$ by detaching the exiting line $\ell_{T}$
and reattaching it to all nodes $w\in N_{2}$, and by detaching the entering
line $\ell'_{T}$ then reattaching it to all nodes $w\in N_{1}$.
One can note that the product of the node factors of a cluster 
$T'\in\calF_{1}(T)$ differs from that of $T$ only because of an extra factor 
$-\nn_{v}\nn_{w}$, where $v\in N_{1}$ is the node which $\ell_{T}$ 
is attached to and $w\in N_{2}$ is the node which $\ell'_{T}$ enters
(again we are considering together all self-energy clusters
with the entering and exiting lines attached to the same nodes, respectively).
Indeed, detaching $\ell_{T}$ from a node $w_{1}\in N_{1}$ and then 
reattaching it to $w_{2}\in N_{1}$, some node factors of the nodes in
$N(\calP(w_{1},w_{2}))$ (we are denoting by $\calP(w_{1},w_{2})$ the path
connecting $w_{1},w_{2}$ and by $N(\calP(w_{1},w_{2}))$ the set of nodes
connected by lines in $\calP(w_{1},w_{2})$)
can change their sign since some lines can change
their direction (see Lemma \ref{lem:4.12}). 
Of course if the components of a line $\ell\in\calP(w_{1},w_{2})$ are inverted,
the corresponding propagator $\calG_{\ell}=
\Psi_{n_{\ell}}(x_{\ell})(x_{\ell}^{2}\uno-\MM^{[n_{\ell}-1]}(x_{\ell}))_{e_{\ell},u_{\ell}}^{-1}$
is replaced by $\Psi_{n_{\ell}}(-
x_{\ell})((-x_{\ell})^{2}\uno-\MM^{[n_{\ell}-1]}
(-x_{\ell}))_{u_{\ell},e_{\ell}}^{-1}$; thanks to Lemma \ref{lem:4.12}, if $e_{\ell}=u_{\ell}=
\be$ or $e_{\ell}=u_{\ell}\in\{\al_{1},\ldots,\al_{d}\}$ the propagator
does not change when considering $\ell$ as a line in $T'$, otherwise it
changes its sign. 
But since one has $u_{\ell_{T}},e_{\ell_{T}'}\in \{\al_1,\ldots,\al_d\}$,
then the number of changes of sign 
(both in the node factors or in the propagators along $\calP({w_{1},w_{2}})$) is
even, so that the overall product does not change sign.

Reasoning as above, we can conclude that the value of a 
a cluster $T''\in\calF_{2}(T)$ differs from that of $T$ only because of a factor
$-\nn_{v}\nn_{w}$, where $v\in N_{1}$ is the node which $\ell'_{T}$ enters
and $w\in N_{2}$ is the node which $\ell_{T}$ exits.

No other changes are produced, except for
the differentiated propagator  which can change sign: the sign changes
for the clusters in $\calF_{1}(T)$ while it remains the same for those
in $\calF_{2}(T)$.  Then by summing over all possible cluster in
$\calF_{1}(T)$ we obtain  $-\nn_{1}\nn_{2}$ 
times a common factor, while by summing over all possible cluster in 
$\calF_{2}(T)$ we obtain $\nn_{1}\nn_2$ times the same common factor, 
so that the overall sum gives zero. Hence (\ref{eq:4.12a}) is proved.

Now pass to (\ref{eq:4.12b}). Given a cluster $T\in\gotR_{n,u,e}$ with 
$e\in\{\al_1,\ldots,\al_d\}$ and $u=\be$ consider all the self-energy
clusters which can be obtained from $T$ by detaching the entering line
$\ell'_{T}$ (note that $e_{\ell_T}=e$) and  reattaching it  to all the nodes 
$v\in N(T)$. Note that again some momenta can change sign, 
but the correponding propagators does not change (again reasoning as done
for the proof of Lemma \ref{lem:4.12} above). 
Hence we obtain a common factor times $\ii\nn_{v}$ where $v$ is the node which 
the exiting line is attached to, so that $\sum_{T}\Val_{T}(0)=0$.

To prove (\ref{eq:4.12c}) one simply notes that it follows
from \eqref{eq:4.12b} and \eqref{eq:4.10c}.

Finally, given a cluster $T\in\gotR_{n,\be,\be}$, consider a contribution to 
$\partial_{x}\Val_{T}(0)$ in which a line $\ell$ is differentiated
(see \eqref{eq:3.13}). The line $\ell$ divides $N(T)$ into two disjoint set
of nodes $N_{1}$ and $N_{2}$ such that $\ell_{T}$ exits a node
$v_{1}\in N_{1}$ and $\ell'_{T}$ enters 
a node $v_{2}\in N_{2}$ i.e.
$N_{2}=\{w\in N(T) : w\prec \ell\}$ and $N_{1}=N(T)\setminus N_{2}$. 
Again, with the same notations as in \eqref{eq:4.13}, one has
$\nn_{1}+\nn_2=0$.  Then consider the cluster obtained by detaching the
exiting line $\ell_{T}$ from $v_{1}$ and reattaching it to the node $v_{2}$, and, 
at the same time, by detaching the entering line $\ell'_{T}$ from $v_{2}$ and 
reattaching to the node $v_{1}$: note that this new cluster again belongs to 
$\gotR_{n,\be,\be}$. Due to this operation, the directions of the line 
along the path connecting $v_{1}$ to $v_{2}$ are reversed, so that for
such lines the momentum $\nn_{\ell}$ is replaced with $-\nn_{\ell}$ 
but the product of the propagators times the node factors does not change.
This means that no overall change
is produced, except for the differentiated propagator which change the sign.
By summing over the two considered cluster we obtain zero because of the
change of sign of the differentiated propagator.
Hence the assertion follows.\EP

\begin{rmk} \label{rmk:4.14}
\emph{
Lemma \ref{eq:4.13} is the counterpart of \eqref{eq:3.17} for the
renormalised self-energies.
}
\end{rmk}

Set $\Theta_{k,\nn,h}^{\RR,n}=\{\theta\in\Theta_{k,\nn,h}^{\RR}:n_{\ell}
\le n\mbox{ for all }\ell\in L(\theta) \}$ and define
\begin{equation}\label{eq:4.14}
\FF^{\RR,n}(\e,\be_{0}):=\sum_{k\ge0}\e^{k}
\sum_{\ttheta\in\Theta_{k+1,{\vzero},\aaa}^{\RR,n}}\VVal(\ttheta),\qquad
G^{\RR,n}(\e,\be_{0}):=\sum_{k\ge0}\e^{k}
\sum_{\theta\in\Theta_{k+1,\vzero,\be}^{\RR,n}}\Val(\theta).
\end{equation}

\begin{lemma}\label{lem:4.15}
Assume $(\e,\be_{0})\in\SSSS_{p}$. Then one has 
$\e\del_{\be_0}G^{\RR,n}(\e,\be_0)=\MM^{[n]}_{\be,\be}(0)+O(\e^{2}{\rm e}^{-C2^{m_{n+1}}})$,
for some positive constant $C$, for all $n\le p$.
\end{lemma}

The proof of the result above essentially follows the lines of the
proof of Lemma 4.12 in \cite{CG2} and Lemma 4.8 in \cite{CG3}.
In particular it does not depend on the Hamiltonian
structure of the equations of motion.

\begin{rmk}\label{rmk:4.16}
\emph{
 From Lemma \ref{lem:4.15} it follows that, if $(\e,\be_{0})\in\SSSS$,
one can define
\begin{equation*}
\MM^{[\io]}(x):=\lim_{n\to\io}\MM^{[n]}(x),\qquad
\GG^{\RR}(\e,\be_0):=\lim_{n\to\io}\GG^{\RR,n}(\e,\be_{0}),
\end{equation*}
with $\GG^{\RR,n}(\e,\be_{0}):=(\FF^{\RR,n}(\e,\be_{0}),G^{\RR,n}(\e,\be_{0}))$
and one has
\begin{equation}\label{eq:4.15}
\MM^{[\io]}_{\be,\be}(0)=\e\partial_{\be_{0}}G^{\RR}(\e,\be_0).
\end{equation}
Note that (\ref{eq:4.15}) is pretty much the same equality provided
by Lemma 4.8 in \cite{CG2}, adapted to the present case.
}
\end{rmk}

\subsection{A suitable assumption: bifurcation equations}\label{sec:4.4}

Here we shall see how to solve the bifurcation equations (\ref{eq:2.2c}) and (\ref{eq:2.2d})
under the assumption that Property 1 is satisfied; again Property 1
assures that all quantities are well defined. We shall see that (\ref{eq:2.2c})
is automatically satisfied, while (\ref{eq:2.2d}) requires for $\be_{0}$ to be properly chosen
as a function of $\e$.

\begin{lemma}\label{lem:4.17}
For any $(\e,\be_{0})\in\SSSS$ one has $\FF^{\RR}(\e,\be_0)=0$.
\end{lemma}

\prova
Consider a tree  $\theta\in\Theta^{\RR}_{k,0,\al_{i}}$
(that is a contribution to $F^{[k-1]}_{i}(\e,\be)$) with root line 
$\ell_{\theta}$ such that $u_{\ell_{\theta}}=\al_i$ (of course $e_{\ell_{\theta}}=
\al_i$), so that the propagator of the root line is 1. Now consider all
trees $\theta'$ obtained $\theta$ by detaching the root line $\ell_{\theta}$
and reattaching it to all nodes $v\in N(\theta)$. 
By detaching $\ell_{\theta}$ from
$v\in N(\theta)$ and reattaching it to another node $w\in N(\theta)$, the
lines $\ell\in \calP(v,w)$ (we are using the same notation as in the
proof of Lemma \ref{lem:4.13}) change their direction. In this case, given a
node $v_{1}\in N(\calP(v,w))\backslash\{v,w\}$, 
call $\ell_{v_{1}},\ell'_{v_{1}}\in\calP(v,w)$
the lines exiting and entering $v_{1}$ respectively. The node factor $\calF_{v_{1}}$ does 
not change its sign if $u_{\ell_{v_{1}}}=e_{\ell'_{v_{1}}}=\be$ or 
$u_{\ell_{v_{1}}},e_{\ell_{v_2}}\in\{\al_{1},\ldots,\al_{d}\}$ when considering $v_{1}$
as a node in $\theta'$, otherwise the sign of $\calF_{v_{1}}$ changes. 
The node factor $\calF_{v}$ does not change its sign only if 
$e_{\ell_{v}}\in\{\al_{1},\ldots,\al_{d}\}$, while the node factor $\calF_{w}$ 
does not change its sign only if $u_{\ell_{w}}\in\{\al_{1},\ldots,\al_{d}\}$. 
Moreover, given a line $\ell\in\calP(v,w)$, thanks to Lemma \ref{lem:4.12} the
corresponding propagator does not change its sign when one considers $\ell$ as a
line of $\theta'$ only if $e_{\ell}=u_{\ell}=\be$ or $e_{\ell},u_{\ell}\in
\{\al_{1},\ldots,\al_{d}\}$. Since one has $u_{\ell_{\theta}}=\al_{i}$ then the
number of changes of sign, of both the propagators
and of the node factors, is even, so that the overall product does not change.
But in this case, the value of $\theta'$ differs from the value of $\theta$
by a factor $\ii\nn_{v}$,
if $v$ is the node which the root line is attached to. The sum of all such
values is zero because $\sum_{v\in N(\theta)}\nn_{v}=0$.

Let us now consider a tree $\theta\in\Theta^{\RR}_{k,0,\al_{i}}$ with 
$u_{\ell_{\theta}}=\al_{j}$ with $j\neq i$ or $u_{\ell_{\theta}}=\be$. In
this case the value of the tree is zero because the propagator of the root
line is $(\uno)_{e_{\ell_{\theta}},u_{\ell_{\theta}}}=0$. Of course we can reason in the
same way for any $i=1,\ldots,d$, therefore the assertion follows.  \EP

Now consider the equation
\begin{equation}\label{eq:4.16}
G^{\RR}(\e,\be_{0})=0.
\end{equation}
One cannot reason as in Lemma \ref{lem:4.17} above, because in principle there
can be nonzero terms since the first order: in such a case, we have to consider \eqref{eq:4.16}
as an implicit function problem and fix $\be_{0}=\be_{0}(\e)$ in a suitable way.

\begin{lemma}\label{lem:4.18}
Assume that there exists $\bar\e >0$ such that $\SSSS=[-\bar\e,\bar\e]\times\TTT$.
Then there exist at least two values $\be_{0}=\be_{0}(\e)$ such that (\ref{eq:4.16}) 
is satisfied for  $\e$ small enough.
\end{lemma}

\prova
Thanks to the variational nature of the Hamilton equations, the function $G^{\RR}$ is the
$\be_{0}$-derivative of the average  of the Lagrangian $\gamma$ computed
along the solution of the range equations (see the comments at the end of
Section \ref{sec:3}).
Under the assumption that Property 1 holds for all $\be_{0}\in\TTT$, 
$\gamma$ is $C^{\infty}$ for any $\be_{0}\in\TTT$ and hence it has at least
two critical points.\EP

If Property $1$ does not hold for all $\be_{0}\in\TTT$
--- or simply if this is not known ---, we have to reason in a different way.
First of all, let us formally expand $G^{\RR}$ in power series in
$\e$, by writing $G^{\RR}(\,\be_{0})=
\sum_{k\geq0}\e^{k}G^{\RR(k)}(\be_{0})$.
Note that $G^{\RR(k)}(\be_{0})$ equals $[\del_{\be}f(\aaa,\be)
]^{(k)}_{\vzero}$ and hence  can be written as a sum over non-renormalised trees
as in (\ref{eq:3.6b}). 

If one has $G^{\RR(k)}(\aaa_{0},\be_{0})\equiv0$ for all $k\geq0$, then
(\ref{eq:4.16}) is formally satisfied.
Otherwise the following condition makes sense.

\begin{cond}\label{cond:2}
Either $G^{\RR(0)}(\be_{0})$ is not identically vanishing or there
exists $k_{0}\in\NNN$ such that $G^{\RR(k)}(\be_{0})\equiv0$ for $0\le k < k_{0}$,
while $G^{\RR(k_{0})}(\be_{0})$ is not identically vanishing.
\end{cond}

\begin{rmk} \label{rmk:4.19}
\emph{
We know that $G^{\RR(k_{0})}$ is the derivative with respect to $\be_{0}$
of the time average of the $k_{0}$-th order Lagrangian $\g^{(k_{0})}$ 
computed along the formal solution.
Since $\g^{(k_{0})}$ is analytic and periodic in $\be_{0}$, 
and it is not identically constant, then it admits at least one
maximum and one minimum.
In particular, for $\s=\pm$, there exist  $\be_{0,\s}^{*}\in\TTT$ and $\gotn_{\s}\in\NNN$,
with $\gotn_{\s}$ odd, such that $(\s1)^{k_{0}+1}\del_{\be_{0}}^{\gotn_{\s}}G^{\RR(k_{0})}(\be_{0,\s}^{*})<0$.
}
\end{rmk}

\begin{rmk} \label{rmk:4.20}
\emph{
Under Condition \ref{cond:2} we can write
\begin{equation} \nonumber 
G^{\RR}(\e,\be_{0})=\e^{k_{0}}\left( G^{\RR(k_{0})}(\be_{0})
+G^{\RR(>k_{0})}(\e,\be_{0})\right),
\end{equation}
where $k_{0} \ge 0$ and $G^{\RR(>k_{0})}(\e,\be_{0})=O(\e)$; hence we can solve the 
equation of motion up to order $k_{0}$ without fixing the parameter $\be_{0}$.
}
\end{rmk}

With the notations in (\ref{eq:3.16}), the condition
that $G^{\RR}(\e,\be_{0})$ identically vanishes to all orders is equivalent
to the condition that $\matL_{\be,\be}^{(k)}\equiv 0$ for all $k\ge 1$
(see comments at the end of Section \ref{sec:3}). Therefore the only
condition left when neither Condition \ref{cond:1} nor Condition \ref{cond:2}
are satisfied is the following.

\begin{cond}\label{cond:3}
One has $\matL^{(k)}_{\be,\be}\equiv0$ for all $k\ge1$ and there exists
$i=1,\ldots,d$ and $k_{1}\in\NNN$ such that $\DD^{(k)}_{\al_{i},\be}\equiv0$
for $k<k_{1}$ while $\DD^{(k_{1})}_{\al_{i},\be}$ does not vanishes
identically.
\end{cond}

\begin{rmk}\label{rmk:4.21}
\emph{
If we take the formal expansion of the functions $\FF^{\RR}
(\e,\be_{0})$, $G^{\RR}(\e,\be_{0})$ and
$\MM^{[\io]}_{u,e}(0)$, $u,e\in\{\al_{1},
\ldots,\al_{d},\be\}$, we obtain the tree expansions of
Section \ref{sec:3}, where the self-energy clusters are allowed.
Then, as we have seen in Lemma \ref{lem:3.5}, the identity (\ref{eq:4.15})
holds to any perturbation order.
If we assume Condition \ref{cond:2} we obtain
\begin{equation}\label{eq:4.17}
\sum_{k=1}^{k_{0}-1}\e^{k} [\MM^{[\io]}_{\be,\be}(0)]^{(k)} \equiv 0
\qquad \Longrightarrow \qquad
\Biggl| \sum_{k=1}^{k_{0}-1}\e^{k} [\MM^{[n]}_{\be,\be}(0)]^{(k)} \Biggr|\le
\e^{2} A_{1} \, {\rm e}^{-A_{2}2^{m_{n}}}, \end{equation}
for some positive constants $A_{1}$ (depending on $k_{0}$) and $ A_{2}$.
If Condition \ref{cond:3} is satisfied, then \eqref{eq:4.17}
is satisfied for any (finite) $k_{0}$, so that
$\MM^{[\io]}_{\be,\be}(0)\to0$ faster than any power as $\e\to0$; moreover in such a case
$[\partial_{x}\MM^{[\io]}_{\aaa,\be}(0)]^{(k)}\equiv0$ 
for all $k=1,\ldots,k_{1}-1$ and
\begin{equation}\label{eq:04.28}
\Biggl| \sum_{k=1}^{k_{1}-1}\e^{k}
[\partial_{x}\MM^{[n]}_{\aaa,\be}(0)]^{(k)}
\Biggr|\le \e^{2}B_{1}\, {\rm e}^{-B_{2}2^{m_{n}}}, 
\end{equation}
for some positive constants $B_{1}$ (depending on $k_{1}$) and $B_{2}$. 
}
\end{rmk}

Assume Condition \ref{cond:2} and fix $\s\in\{\pm1\}$.
Suppose for the time being $\SSSS$ to be
an open set containing $(0,\be_{0,\s}^{*})$.
Then, by reasoning as for Lemma 4.15 of \cite{CG2}, one can show that
(i) there exists a neighbourhood $U$ of $\e=0$ such that the implicit function equation (\ref{eq:4.16}) admits in $\SSSS$
a solution $\be_{0}=\be_{0,\s}(\e)$, with $\e\in U$ and $\be_{0,\s}(0)=\be_{0,\s}^{*}$; (ii) for $\hbox{sign}\,\e=\s1$
one has $\e\del_{\be_{0}}G^{\RR}(\e,\be_{0,\s}(\e))\leq0$.
Then for $\e\in U$, with $\hbox{sign}\,\e=\s1$, and $\be_{0}=\be_{0,\s}(\e)$,
the functions $\aaa^{\RR}$, $\be^{\RR}$
in \eqref{eq:4.7} are well defined and one has
\begin{equation*}\begin{aligned}
& \FF^{\RR}(\e,\be_0)=
[-\del_{\aaa}f(\aaa^{\RR}(t;\e,\aaa_0,\be_0),
\be^{\RR}(t;\e,\aaa_0,\be_0))]_{{\vzero}},\\
& G^{\RR}(\e,\be_0)=
[\del_{\be}f(\aaa^{\RR}(t;\e,\aaa_0,\be_0),
\be^{\RR}(t;\e,\aaa_0,\be_0))]_{\vzero} ,
\end{aligned}
\end{equation*}
and hence by Lemma \ref{lem:4.10} the functions $\aaa^{\RR}(t;\e,\aaa_0,\be_{0,\s}(\e))$
and $\be^{\RR}(t;\e,\aaa_{0},\be_{0,\s}(\e))$ solve the equation of motion ($\ref{eq:1.2}$).
However, the argument above is not sufficient to prove 
the existence of a quasi-periodic solution with frequency $\oo$,
because we have assumed --- without proving --- 
that Property 1 is satisfied on a non-empty open set. In Section \ref{sec:4.5} we shall show that, 
thanks to the simmetry property of Lemma $\ref{lem:4.10}$ and 
the identity of Lemma $\ref{lem:4.15}$, Property 1 is
satisfied along a suitable curve 
$\be_{0}=\ol{\be}_{0}(\e)$ such that $G^{\RR}(\e,\ol{\be}_{0}(\e))=0$ and
$\ol{\be}_{0}(\e)$ is continuous for $\e\neq0$.
More precisely, we shall proceed by induction as follows. Under Condition \ref{cond:3},
assuming that Property 1-$n$ holds for all $n<p$ will imply, thanks to the bounds and symmetry properties
seen in the previous sections, that also Property 1-$p$ holds. The discussion of Condition \ref{cond:2}
is more delicate: we shall need to introduce some auxiliary quantities for which an analogous result
is obtained and then show that this yields the same result for the self-energies.

\subsection{Convergence of the resummed series}
\label{sec:4.5}

First of all we recall that if we formally expand the resummed series,
we obtain the same formal expansion
as in Section \ref{sec:3}. In particular, either Condition \ref{cond:1}
is satisfied --- and hence we can reason as in Section \ref{sec:3} --- 
or at least one among $\matL^{(k)}_{\be,\be}$
and $\DD^{(k)}_{\al_{i},\be}$ for $i=1,\ldots,d$ is not identically vanishing.
Let us start from the case in which Condition \ref{cond:3} holds.

\begin{lemma}\label{lem:4.22}
Assume Condition \ref{cond:3}.  Then $\MM$ satisfies Property $1$ for all $\be_{0}\in\TTT$
and $\e$ small enough.
\end{lemma}

\prova
We shall prove that $\MM$ satisfies Property 1-$p$ 
for all $p\geq0$, by induction on $p$.
Property 1-$0$ is trivially satisfied for $\e$ small enough.
Indeed the matrix
$\MM^{[-1]}(x)$ defined in \eqref{eq:04.16}
is the null matrix, so that $\calG^{[0]}(x)=\uno\Psi_{0}(x)/x^{2}$,
and hence $\Arrowvert\calG^{[0]}(x)\Arrowvert \leq c_{0}/x^{2(d+1)}$,
for some constant $c_{0}>0$.
Assume that $\MM$ satisfies Property 1-$p$. By Lemmas 
\ref{lem:4.8} and \ref{lem:4.13}
\begin{equation} \nonumber 
\MM^{[p]}(x)=\left(\begin{matrix}0_{d} & \vzero \\ 
\vzero & \MM^{[p]}_{\be,\be}(0)\end{matrix}\right)+
x\left(\begin{matrix}0_{d} & \del_{x}\MM^{[p]}_{\aaa,\be}(0) \\ 
\del_{x}\MM^{[p]}_{\be,\aaa}(0)& 0\end{matrix}\right)+
O(\e^{2}x^{2}).
\end{equation}
We have to bound from below the determinant of the matrix 
$x^{2}\uno-\MM^{[p]}(x)$: we have
\begin{equation} \label{eq:4.19}
\det(x^{2}\uno-\MM^{[p]}(x))=
x^{2d}\Bigl(x^{2}-\Bigl(\MM^{[p]}_{\be,\be}(0)-
\left|\del_{x}\MM^{[p]}_{\aaa,\be}(0)\right|^{2}\Bigr)
+O(\e^{2}x^{2})\Bigr),
\end{equation}
so that we have to show that
\begin{equation} \nonumber 
\MM^{[p]}_{\be,\be}(0)\leq\frac{x^{2}}{2}+
\left|\partial_{x}\MM_{\aaa,\be}^{[p]}(0)\right|^{2},
\qquad  \hbox{with }\quad x^{2} \ge \frac{\al_{m_{p+1}}(\oo)^{2}}{2^{8}}.
\end{equation}
Since
\begin{equation} \nonumber
\left| \partial_{x}\MM_{\aaa,\be}^{[p]}(0) \right| \le
\left|
\sum_{k=1}^{k_{1}} \e^{k} [\partial_{x}\MM_{\aaa,\be}^{[p]}(0)]^{(k)}
\right| +O(\e^{k_{1}+1})
\leq \e^{2} B_{1}{\rm e}^{-B_{2}2^{m_{p}}}+O(\e^{k_{1}+1}),
\end{equation}
and (use Remark \ref{rmk:4.21} with $k_{0}=2k_{1}+2$)
\begin{equation}\nonumber
\left| \MM_{\be,\be}^{[p]}(0) \right|
\le\left|\sum_{k=1}^{2k_{1}+2} \e^{k} [\MM_{\be,\be}^{[p]}(0)]^{(k)}\right|
+O(\e^{2k_{1}+3})
\leq \e^{2} A_{1} \, {\rm e}^{-A_{2}2^{m_{p}}}+O(\e^{2k_{1}+3}),
\end{equation}
the assertion follows by the condition $\BB(\oo)<\infty$.\EP

By Lemma \ref{lem:4.22} we can apply Lemma \ref{lem:4.18} and deduce,
in the case of Condition \ref{cond:3}, the existence of
at least two $d$-dimensional invariant tori.
Therefore  we are left with Condition \ref{cond:2}.

First of all, for all $n\geq0$, we define the $C^{\infty}$ non-decreasing
functions $\x_{n}$ such that
\begin{equation}\label{eq:4.20}
\x_{n}(x):=\left\{
\begin{aligned}
&1, \quad x\leq\al_{m_{n+1}}(\oo)^{2}/2^{12},\\
&0,\quad x\geq\al_{m_{n+1}}(\oo)^{2}/2^{11},
\end{aligned}\right.
\end{equation}
and set $\x_{-1}(x)=1$. Define recursively, for all $n\geq0$, the
\emph{regularised propagators}
\begin{equation} \nonumber 
\ol\calG^{[n]}(x):=\Psi_{n}(x)\left(x^{2}\uno-
\ol\MM^{[n-1]}(x)\x_{n-1}(\Delta_{n-1})\right)^{-1}
\end{equation}
with $\ol\MM^{[-1]}(x)=\MM^{[-1]}(x)$ as given by \eqref{eq:04.16}
and, for all $n\geq0$,
\begin{equation} \nonumber 
\ol\MM^{[n]}(x):=\ol\MM^{[n-1]}(x)+
\chi_{n}(x)\ol{M}^{[n]}(x),
\end{equation}
where we have set for all $u,e\in\{\al_{1},\ldots,\al_{d},\be\}$,
$$
\ol{M}^{[n]}_{u,e}(x):=
\sum_{T\in\gotR_{n,u,e}}\e^{k(T)}\TVal_{T}(x) ,
$$
with
$$
\TVal_{T}(x):=\left(\prod_{v\in N(T)}\calF_{v}\right)
\left(\prod_{\ell\in L(T)}\ol\calG^{[n_{\ell}]}_{e_{\ell},u_{\ell}}(\oo\cdot\nn_{\ell})\right)
$$
and
\begin{equation} \nonumber 
\Delta_{n-1}=\Delta_{n-1}(\e,\be_{0}):=
\ol\MM^{[n-1]}_{\be,\be}(0;\e,\be_{0})-
\sum_{k=0}^{k_{0}-1}\e^{k}[\ol\MM^{[n-1]}_{\be,\be}(0;\e,
\be_{0})]^{(k)}.
\end{equation}
Set also $\ol\MM:=\{\ol\MM^{[n]}(x)\}_{n\geq-1}$ and 
$\ol\MM^{\x}:=\{\ol\MM^{[n]}(x)\xi_{n}(\Delta_{n-1})\}_{n\geq-1}$.

\begin{lemma}\label{lem:4.23} 
$\ol\MM^{\x}$ satisfies Property 1 for $\e$ small enough and any $\be_{0}\in\TTT$.
\end{lemma}

\prova
We shall prove that $\ol\MM^{\x}$ satisfies Property 1-$p$ 
for all $p\geq0$, by induction on $p$.
Property 1-$0$ is trivially satisfied for $\e$ small enough.
Indeed the matrix $\ol\MM^{[-1]}(x)$
is self-adjoint, so that also $\ol\calG^{[0]}(x)$ is
self-adjoint and we can estimate its eigenvalues and conclude 
$\Arrowvert\ol\calG^{[0]}(x)\Arrowvert
\leq \bar{c}_{0}/x^{2(d+1)}$, for some $\bar{c}_{0}>0$.
Assume then that $\ol\MM^{\x}$ satisfies Property 1-$p$. Then we can repeat
almost word by word the proof of Lemmas \ref{lem:4.8} and \ref{lem:4.13},
as done in \cite{CG2,CG3} so as to obtain
\begin{equation}\label{eq:4.21}
\ol\MM^{[p]}(x)=\left(\begin{matrix}0_{d} & \vzero \\ 
\vzero & \ol\MM^{[p]}_{\be,\be}(0)\end{matrix}\right)+
x\left(\begin{matrix}0_{d} & \del_{x}\ol\MM^{[p]}_{\aaa,\be}(0) \\ 
\del_{x}\ol\MM^{[p]}_{\be,\aaa}(0)& 0\end{matrix}\right)+
O(\e^{2}x^{2}).
\end{equation}
We have to bound from below the determinant of the matrix 
$x^{2}\uno-\ol\MM^{[p]}(x)\x_{p}(\Delta_{p})$.
 From \eqref{eq:4.21} it is easy to check that such determinant is
\begin{equation} \label{eq:4.22}
x^{2d}\Bigl(x^{2}-\Bigl(\ol\MM^{[p]}_{\be,\be}(0)-
\left|\del_{x}\ol\MM^{[p]}_{\aaa,\be}(0)\right|^{2}\Bigr)
\x_{p}(\Delta_{p}) + O(\e^{2}x^{2})\Bigr).
\end{equation}
Thanks to the definition of the functions $\xi_{p}$, since
$$
\sum_{k=0}^{k_{0}-1}\e^{k}[\ol\MM^{[p]}_{\be,\be}(0)]^{(k)}
= \sum_{k=0}^{k_{0}-1} \e^{k} [\MM^{[p]}_{\be,\be}(0)]^{(k)}$$
by Remark \ref{rmk:4.21}, one has
\begin{equation*}
x^{2}-\Bigl( \ol\MM^{[p]}_{\be,\be}(0)-
\left| \del_{x}\ol\MM^{[p]}_{\aaa,\be}(0)\right|^{2} \Bigr)\x_{p}
(\Delta_{p})\geq
x^{2}-\Bigl(\ol\MM^{[p]}_{\be,\be}(0)\Bigr)\x_{p}
(\Delta_{p})\ge \frac{x^{2}}{2} .
\end{equation*}
Then $\|\ol\calG^{[p+1]}(x)\|\leq \bar{c}_{1}/x^{2(d+1)}$,
for some positive constant $\bar{c}_{1}$, that is Property 1-$(p+1)$ with
$c_{2}=2(d+1)$ in Definition \ref{def:4.4}.\EP

Set
\begin{equation}\label{eq:4.23}
\ol\aaaa_{\nn}^{[k]}(\e;\aaa_{0},\be_{0}) :=
\sum_{\ttheta\in\Theta^{\RR}_{k,\nn,\aaa}}
\TTVal(\ttheta), \qquad
\ol{b}_{\nn}^{[k]}(\e;\aaa_{0},\be_{0}) := \sum_{\theta\in\Theta^{\RR}_{k,\nn,\be}}
\TVal(\theta),\qquad
\nn\ne\vzero,
\end{equation}
where, for $\ttheta=(\theta_{1},\ldots,\theta_{d})\in
\Theta^{\RR}_{k,\nn,\aaa}$ we denoted $\TTVal(\ttheta):=
(\TVal(\theta_{1}),\ldots,\TVal(\theta_{d}))$, and define
\begin{equation}\label{eq:4.24}
\begin{aligned}
&\ol\aaaa(t;\e,\aaa_{0},\be_{0})=\sum_{k\geq1}\e^{k}
\sum_{\nn\in\ZZZ^{d}_{*}}
{\rm e}^{\ii\nn\cdot\oo t}\ol\aaaa_{\nn}^{[k]}, \qquad
\ol{b}(t;\e,\aaa_{0},\be_{0})=\sum_{k\geq1}\e^{k}\sum_{\nn\in\ZZZ^{d}_{*}}
{\rm e}^{\ii\nn\cdot\oo t}\ol{b}_{\nn}^{[k]}, \\
&\ol{G}(\e,\be_{0}):=\sum_{k\ge0}\e^{k}\ol{G}^{(k)}(\e,\be_{0}):=
\sum_{k\ge0}\e^{k}\sum_{\theta\in\Theta^{\RR}_{k+1,\vzero,\be}}
\TVal(\theta).
\end{aligned}
\end{equation}
A result analogous to Lemma \ref{lem:4.9} holds and can be proved
in the same way (see \cite{CG2,CG3}), so we conclude that
the series \eqref{eq:4.24} converge. However, because of the presence of
the functions $\xi_{n}$, in principle no equivalent of Lemma \ref{lem:4.10}
applies in this case. In other words, in general  the functions (\ref{eq:4.24})
are no longer solutions of the equations of motions, unless
$\xi_{n}(\Delta_{n})\equiv1$.  Therefore we would like
to show that, for any $\e$ small enough, it is possible to fix suitably
 $\be_{0}=\ol{\be}_{0}(\e)$ in such a way that $\xi_{n}(\Delta_{n})$ be
 identically one.
 
\begin{lemma}\label{lem:4.24}
One has
$[\ol{G}(\e,\be_{0})]^{(k)}=[G^{\RR}(\e,\be_{0})]^{(k)}$
for all $k=0,\ldots,k_{0}$.
\end{lemma}

\proof
Set $\Theta_{k,\nn,\be}^{\RR(n)}:=\{\theta\in\Theta_{k,\nn,\be}^{\RR,n} :
\exists \ell\in L(\theta) \hbox{ such that }n_{\ell}=n\}$ and write
\begin{equation}\nonumber
\ol{G}(\e,\be_{0})=\sum_{k\ge 0}\e^{k}\sum_{n\ge0}
\sum_{\theta\in\Theta^{\RR(n)}_{k+1,\vzero,\be}}\TVal
(\theta) .
\end{equation}
Note that if $\theta\in\Theta_{k,\vzero,\be}^{\RR(n)}$ one has
$\prod_{v\in N(\theta)}|\calF_{v}|\le E_{1}^{k} 
{\rm e}^{-E_{2}2^{m_{n}}}$,
for some constants $E_{1},E_{2}$. Moreover one can write formally
$$
\ol{\calG}^{[n_{\ell}]}(x) = \Psi_{n_{\ell}}(x) \frac{1}{x^{2}} \Bigl( \uno +
\sum_{m\ge1}\Bigl( \frac{1}{x^{2}} \ol{\MM}^{[n_{\ell}-1]}
(x)\x_{n_{\ell}-1}(\Delta_{n_{\ell}-1}) \Bigr)^{m}\Bigr),
$$
and
$\x_{n_{\ell}-1}(\Delta_{n_{\ell}-1}) = 1+\x_{n_{\ell}-1}'(\Delta^{*})
\Delta_{n_{\ell}-1}$ for some $\Delta^{*}$, where $\Delta_{n_{\ell}-1}=O(\e^{k_{0}})$
and
$$
|\x_{n_{\ell}-1}'(\Delta^{*})|\le \frac{E_3}{\al_{m_{n_{\ell}}}(\oo)^2}\le 
\frac{E_3}{\al_{m_{n}}(\oo)^2},
$$
for some positive constant $E_{3}$ independent of $n$. Hence the assertion
follows.
\EP

Define
\begin{equation}\label{eq:4.25}
\ol\MM^{[\infty]}(x):=\lim_{n\to\infty}
\ol\MM^{[n]}(x),
\end{equation}
and note that, by Lemma \ref{lem:4.23}, the limit in \eqref{eq:4.25}
is well defined, and it is $C^{\infty}$ in both $\e$ and $\be_{0}$.

For $\s=\pm$ let us  introduce the $C^{\infty}$
functions $R(\e,\be_{0})$ such that 
$\ol\MM^{[\infty]}_{\be,\be}(0)=\e\del_{\be_0}R(\e,\be_{0})$.
Note that $R(\e,\be_{0})=\e^{k_{0}}\Gamma(\e,\be_{0})$, with
$\Gamma(\e,\be_{0})=G^{\RR(k_{0})} (\be_{0})+O(\e)$,
so that $\Gamma(0,\be_{0,\s}^{*})=0$ and
$\del_{\be_{0}}^{\gotn_{\s}}\Gamma_{\s}(0,\be_{0,\s}^{*})\neq0$,
with $\be_{0,\s}^{*}$ and $\gotn_{\s}$ defined in Remark \ref{rmk:4.19}.
For any of such function consider the implicit function equation
\begin{equation}\label{eq:4.26}
R(\e,\be_0)=0.
\end{equation}
%
\begin{lemma}\label{lem:4.25}
Assume Condition \ref{cond:2}.
There exist a neighbourhood $U$ of $\e=0$ and, for $\e\in U$,
a solution $\be_{0}=\ol\be_{0}(\e)$ to the implicit function equation (\ref{eq:4.26}), such that
$$ \lim_{\e\to 0^{\s}} \ol\be_{0}(\e)=\be_{0,\s}^{*} , \quad \s=\pm , \qquad
\e\del_{\be_{0}}R(\e,\ol{\be}_{0}(\e))\leq0  . $$ 
Moreover $\ol{\be}_{0}(\e)$ is continuous in $U$ for $k_{0}$ odd 
and in $U\setminus\{0\}$ for $k_{0}$ even.
\end{lemma}

\prova By construction, all the functions $\Gamma(\e,\be_{0})$ are smooth
for $\be_{0}\in \TTT$ and $\e$ small enough. 
Then there exist two half-neighbourhoods $V_{\s,-}$ and $V_{\s,+}$ of $\be_{0}=\be_{0,\s}^{*}$ such that
$\Gamma(0,\be_{0})>0$ for $\be_{0}\in V_{\s,+}$ and 
$\Gamma(0,\be_{0})<0$ for $\be_{0}\in V_{\s,-}$. 
By continuity, there exist
a neighbourhood $U_{\s}=(-\bar\e_{\s},\bar\e_{\s})$ and a continuous
curve $\be_{0,\s}(\e)$ such that
$\be_{0,\s}(0)=\be_{0,\s}^{*}$ and $\Gamma(\e,\be_{0,\s}(\e))\equiv 0$ for $\e\in U_{\s}$.
Moreover if $\del_{\be_{0}}^{\gotn_{\s}}G^{\RR(k_{0})}(\be_{0,\s}^{*})>0$, 
then $V_{\s,+}$ and $V_{\s,-}$ are of the form $(\be_{0,\s}^{*}, v_{\s,+})$ and $(v_{\s,-}, \be_{0,\s}^{*})$,
respectively, and therefore $\del_{\be_{0}}\Gamma(\e,\be_{0,\s}(\e))\geq0$
for all $\e\in U_{\s}$. If on the contrary $\del_{\be_{0}}^{\gotn_{\s}}G^{\RR(k_{0})}(\be_{0,\s}^{*})<0$,
one has $V_{\s,+}=(v_{+},\be_{0,\s}^{*})$ and $V_{\s,-}=(\be_{0,\s}^{*},v_{-})$, and then
$\del_{\be_{0}}\Gamma(\e,\be_{0,\s}(\e))\leq0$ for all $\e\in U_{\s}$. 

If $k_{0}$ is odd, then $\be_{0,+}^{*}=\be_{0,-}^{*}$ and hence one can take $\ol{\be}_{0}(\e)=
\be_{0,+}(\e)=\be_{0,-}(\e)$ in such a way that it is a continuous function of $\e\in U_{+}=U_{-}$.
If $k_{0}$ is even, then one has $\ol{\be}_{0}(\e)=\be_{0,+}(\e)$ for $\e>0$ and
$\ol{\be}_{0}(\e)=\be_{0,-}(\e)$ for $\e<0$, so that $\ol{\be}_{0}(\e)$ has a discontinuity at $\e=0$.\EP

\begin{lemma}\label{lem:4.26}
Assume Condition \ref{cond:2}.
Let $U$ and $\ol\be_{0}(\e)$ be the neighbourhood and the solution
referred to in Lemma \ref{lem:4.25}, respectively.
Then whenever $\e\in U$ and $\be_{0}=\ol{\be}_{0}(\e)$, one has $\x_{n}(\Delta_{n})\equiv1$ 
for all $n\geq0$.
\end{lemma}

\prova By Lemma \ref{lem:4.25},
for $\e\in U$ and $\be_{0}=\ol\be_{0}(\e)$, one has 
$\ol\MM^{[\infty]}_{\be,\be}(0) =
\e\del_{\be_{0}}R(\e,\ol\be_{0}(\e))$. Hence, since the matrices
$\ol{M}^{[n]}(x)$ satisfy bounds analogous to those in Lemma \ref{lem:4.6},
possibly renaming the constants, one has for $\be=\ol{\be}_{0}(\e)$
\begin{equation} \nonumber
\begin{aligned}
& \ol\MM^{[n]}_{\be,\be}(0)-
\sum_{k=1}^{k_{0}-1}\e^{k}[\ol\MM^{[n]}_{\be,\be}(0)]^{(k)}
\leq\ol\MM^{[n]}_{\be,\be}(0)-
\ol\MM^{[\infty]}_{\be,\be}(0)+ \e^{2} \, A_{1} {\rm e}^{-A_{2}2^{m_{n+1}}}
\\
& \quad \leq\sum_{p\geq n+1}|\ol{M}^{[p]}_{\be,\be}(0)|
+\e^{2} A_{1}{\rm e}^{-A_{2}2^{m_{n+1}}}
\leq2K_{0}\e^{2} {\rm e}^{-\ol{K}_{0}2^{m_{n+1}}}+\e^{2} A_{1}{\rm e}^{-A_{2}2^{m_{n+1}}}
\leq\frac{\al^{2}_{m_{n+1}}(\oo)}{2^{13}},
\end{aligned}
\end{equation}
so the assertion follows by the definition of $\x_{n}$.\EP

The following result concludes the proof of the existence of an invariant $d$-dimensional torus
under Condition \ref{cond:2}.

\begin{lemma}\label{lem:4.27}
Assume Condition \ref{cond:2} and
let $\ol{\be}_{0}(\e)$ be as in Lemma \ref{lem:4.25}.
One can choose
the function $R(\e,\be_{0})$ such that
$R(\e,\ol{\be}_{0}(\e))=G^{\RR}(\e,\ol{\be}_{0}(\e))\equiv0$, where
\begin{equation}\nonumber
G^{\RR}(\e,\ol{\be}_{0}(\e)):=
\lim_{n\to\io}G^{\RR,n}(\e,\ol{\be}_{0}(\e))
\end{equation}
and the functions $G^{\RR,n}$ are defined in (\ref{eq:4.14}).
 In particular $(\aaa(t,\e),\be(t,\e))=(\aaa_{0}+\oo t,\ol{\be}_{0}(\e))+
(\aaaa^{\RR}(t;\e,\aaa_{0},\ol{\be}_{0}(\e)),b^{\RR}(t;\e,\aaa_{0},
\ol{\be}_{0}(\e)))$ defined in (\ref{eq:4.7}) solves
the equation of motion (\ref{eq:1.2})
\end{lemma}

\prova
It follows from the results above. Indeed, for any primitive $R$
there is a curve $\ol{\be}_{0}(\e)$ along which
$\MM=\ol{\MM}=\ol{\MM}^{\xi}$ (hence
$\MM$ satisfies Property $1$) and $R(\e;\ol{\be}_{0}(\e))\equiv 0$.
By Lemma \ref{lem:4.15} and the fact that $\MM$ satisfies
Property $1$, also $G^{\RR}$ is among the primitives of
$\MM^{[\io]}_{\be,\be}$ and hence the assertion follows.\EP

\begin{rmk} \label{rmk:4.28}
\emph{
If one considers a convex unperturbed Hamiltonian, e.g. with a plus sign instead of
the minus sign in (\ref{eq:1.1}), one can try to proceed in the same way.
Some parts of the construction simplify: for instance, the
self-energies $\MM^{(k)}(x,n)$ turn out to be self-adjoint and
$(\MM^{(k)}(x,n))^{T}=\MM^{(k)}(-x,n)$. On the other hand,
when dealing with Conditions \ref{cond:2} and \ref{cond:3},
one has to bound from below determinants which have the form
(\ref{eq:4.22}) or (\ref{eq:4.19}), respectively, with the major
difference that a sign plus appears in front of the
squared term; for instance (\ref{eq:4.19}) becomes
$$ x^{2d}\Bigl(x^{2}-\Bigl(\MM^{[p]}_{\be,\be}(0) +
\left|\del_{x}\MM^{[p]}_{\aaa,\be}(0)\right|^{2}\Bigr)
+O(\e^{2}x^{2})\Bigr). $$
Then information on the sign of $\MM^{[p]}_{\be,\be}(0)$
is not enough to control the corrections to $x^{2}$ and hence no lower bound
follows for the determinant. 
Therefore in order to recover Cheng's result further cancellations seem
to be necessary. In turn this means that one should expect
other symmetries to hold for the self-energies.
}
\end{rmk}

\noindent{\bf Acknowledgments.} We are indebted to Pavel Plotnikov
for useful discussions.


\end{document}